\documentclass[9pt,11pt]{article}%
\usepackage{amsfonts}
\usepackage[applemac]{inputenc}
\usepackage{amssymb}
\usepackage{makeidx}
\usepackage{graphicx}
\usepackage{amsmath}
\usepackage[a4paper]{geometry}
\usepackage[displaymath]{lineno}
\usepackage[singlespacing]{setspace}
\usepackage[normalem]{ulem}
\usepackage{lineno}
\usepackage{color}%
\providecommand{\U}[1]{\protect\rule{.1in}{.1in}}
\providecommand{\U}[1]{\protect\rule{.1in}{.1in}}

\newtheorem{theo}{Theorem}
\newtheorem{lem}[theo]{Lemma}

\newtheorem{cor}[theo]{Corollary}
\newtheorem{rem}{Remark}
\newtheorem{exam}{Example}

\newenvironment{dem}[1][Proof]{\noindent \textbf{#1.} }{\ \rule{0.5em}{0.5em}}
\begin{document}

\title{Valadier-like formulas for the supremum function I\thanks{Research of the
first and the second authors is supported by CONICYT grants, Fondecyt no.
1150909 and 1151003, Basal PFB-03 and Basal FB003. Research of the first and
third authors supported by MINECO of Spain and FEDER of EU, grant
MTM2014-59179-C2-1-P. Research of the third author is also partially supported
by the Australian Research Council: Project DP160100854.}}
\author{R. Correa\thanks{e-mail: rcorrea@dim.uchile.cl}, A.\ Hantoute\thanks{e-mail:
ahantoute@dim.uchile.cl (corresponding author)}, M.A. López\thanks{e-mail:
marco.antonio@ua.es}\\${^{\dag}}${\scriptsize Universidad de O'Higgins, Chile, and DIM-CMM of
Universidad de Chile}\\$^{\ddag}${\scriptsize Center for Mathematical Modeling (CMM), Universidad de
Chile }\\$^{§}${\scriptsize Universidad de Alicante, Spain, and CIAO, Federation
University, Ballarat, Australia}}
\date{}
\maketitle

\begin{abstract}
We generalize and improve the original characterization given by Valadier
\cite[Theorem 1]{Valadier} of the subdifferential of the pointwise supremum of
convex functions, involving the subdifferentials of the data functions at
nearby points. We remove the continuity assumption made in that work and
obtain a general formula for\ such a subdifferential. In particular, when the
supremum is continuous at some point of its domain, but not necessarily at the
reference point, we get a simpler\ version which gives rise to the Valadier
formula. Our starting result is the characterization\ given in\ \cite[Theorem
4]{HanLopZal2008}, which uses the $\varepsilon$-subdifferential at the
reference point.

\textbf{Key words. }Pointwise supremum function, convex functions, Fenchel
subdifferential, Valadier-like formulas.

\emph{Mathematics Subject Classification (2010)}: 26B05,\emph{\ }26J25, 49H05.

\end{abstract}

\section{Introduction\label{Sect1}}

Consider\ a family of extended real-valued convex functions$\ f_{t}%
:X\rightarrow\overline{\mathbb{R}}:=\mathbb{R}\cup\{+\infty\},$ indexed by an
arbitrary\ index set $T,$ and defined in a locally convex topological vector
space $X.$ Under the continuity of the supremum function
\begin{equation}
f:=\sup_{t\in T}f_{t} \label{suprefunction}%
\end{equation}
at a point $x\in X,$ a remarkable pioneering\textbf{\ }result due to Valadier
\cite{Valadier} states that the subdifferential of $f$ at $x$ is completely
characterized by means of the subdifferentials of data functions $f_{t}$ at
nearby points; more precisely
\begin{equation}
\partial f(x)=\bigcap_{\varepsilon>0,\text{ }p\in\mathcal{P}}\overline
{\operatorname*{co}}\left\{  \bigcup\limits_{t\in T_{\varepsilon}(x),\text{
}p(y-x)\leq\varepsilon}\partial f_{t}(y)\right\}  ,
\label{Valadier-like formula}%
\end{equation}
where $\mathcal{P}$ is the family of continuous seminorms in $X,$ and
\begin{equation}
T_{\varepsilon}(x):=\{t\in T\mid f_{t}(x)\geq f(x)-\varepsilon\};
\label{T-epsilon}%
\end{equation}
when $f(x)=+\infty,$ we have $f_{t}(x)=+\infty$ for all $t\in T_{\varepsilon
}(x).$ We also use the notation
\begin{equation}
T(x):=\{t\in T\mid f_{t}(x)=f(x)\}. \label{Tdex}%
\end{equation}

There have been many contributions\ to this problem, in various settings,
depending on the structure of the space $X,$ the algebraic and topological
properties of the index set $T$, the behavior of the function $f(x,t):=f_{t}%
(x),$ etc. See, for instance, \cite{Brondsted72, booktikhomirov2,
booktikhomirov, HanLop08, Danskin67, Ioffe72, LiNg11, Soloviev01,
ZalinescuBook}, among many others.

More recently, using the machinery of $\varepsilon$-subdifferentials, and
under the following condition
\begin{equation}
\operatorname*{cl}f=\sup_{t\in T}\operatorname*{cl}f_{t}%
,\label{closedness criterion}%
\end{equation}
involving the lower semicontinuous (lsc, for short) envelopes, it has been
established in \cite[Theorem 4]{HanLopZal2008} (see also \cite{Han06}) that,
for every $x\in X,$
\begin{equation}
\partial f(x)=\bigcap_{\substack{\varepsilon>0\\L\in\mathcal{F}(x)}%
}\overline{\operatorname*{co}}\left\{  \bigcup\limits_{t\in T_{\varepsilon
}(x)}\partial_{\varepsilon}f_{t}(x)+\mathrm{N}_{L\cap\operatorname*{dom}%
f}(x)\right\}  ,\label{formula using epsilon-subdifferential}%
\end{equation}
where
\begin{equation}
\mathcal{F}(x):=\{L\subset X\mid\ L\text{ is a finite-dimensional liner
subspace such that }x\in L\}.\label{Fdex}%
\end{equation}
When $T$ is finite and $T(x)=T,$ the last formula gives rise to
(\cite{Brondsted72}; see, also, \cite[Corollary 12]{HanLopZal2008})
\begin{equation}
\partial f(x)=\bigcap_{\varepsilon>0}\overline{\operatorname*{co}}\left\{
\bigcup\limits_{t\in T(x)}\partial_{\varepsilon}f_{t}(x)\right\}
,\label{marco5}%
\end{equation}
being the first result of this kind that does not require any continuity conditions.

Both approaches are based on enlargements of $\partial f_{t},$ using either
the exact subdifferentials at nearby points or the $\varepsilon$%
-subdifferentials at the reference point. They are comparable because: (i)
both formulas coincide when the $f_{t}$'s are affine\textbf{,} (ii) both
enlargements are related in the Banach spaces setting thanks to the
Brøndsted-Rockafellar theorem.

The main purpose of the present paper consists of generalizing and improving
formula (\ref{Valadier-like formula}). We drop\ the continuity assumption and
obtain general formulas for\ $\partial f(x),$ in the settings of Banach spaces
(\ref{FVB}) and locally convex spaces (\ref{marco1}), both of them involving
enlargements of $\partial f_{t}$, which use the exact subdifferentials of the
$f_{t}$'s at nearby points.

The paper is organized as follows. After Section 2, devoted to notation and
preliminaries, Section 3 provides the free-continuity extension of
(\ref{Valadier-like formula}) in a Banach space (\ref{FVB}), via the
application of Brøndsted-Rockafellar theorem. Our characterization is
comparable with the main limiting-like results given in \cite{Olopez}. The
validity limitation of this result\ outside the Banach setting is illustrated
in this section by means of some examples. In Section 4 we establish the
second free-continuity formula for $\partial f$ (\ref{marco1}) in the setting
of locally convex spaces. Section 5 is devoted to deriving formula
(\ref{marco2}), which constitutes a simpler version of\ (\ref{FVB}) and
(\ref{marco1}), which is valid only when we are confined to the case in which
the supremum function $f$ is finite and continuous at some point. This last
section finishes with the derivation of (\ref{Valadier-like formula}) from our
formula (\ref{marco2}).

In a forthcoming paper we establish the corresponding counterparts of the
formulas provided in the current paper to the case in which the index set $T$
is compact and the functions\ $t\rightarrow f_{t}(z),$ $z\in
\operatorname*{dom}f,$ are upper semicontinuous on $T_{\varepsilon}(x).$ This
will generalize the second theorem in \cite{Valadier}. In the same paper, some
applications to the characterization of normal cones to sublevel sets are also
given, extending some recent results in \cite{AsAh16}.\textbf{ }

\section{Notation and preliminaries\label{Sect2}}

In this paper $X$ stands for a (real) separated locally convex (lcs, shortly)
space\textbf{,} whose\textbf{ }topological dual space is denoted by $X^{\ast}$
and, unless otherwise stated, is endowed with the weak*-topology.$\ $Hence,
the $(X,X^{\ast})$ forms a dual topological pair by means of the canonical
bilinear form $\langle x,x^{\ast}\rangle=\langle x^{\ast},x\rangle:=x^{\ast
}(x),$ $(x,x^{\ast})\in X\times X^{\ast}.$ The zero vectors in all\textbf{
}the involved spaces are denoted by $\theta$, and the convex, closed and
balanced neighborhoods of $\theta$ are called $\theta$-neighborhoods. The
families of such $\theta$-neighborhoods in $X$ and in $X^{\ast}$ are denoted
by $\mathcal{N}_{X}$ and $\mathcal{N}_{X^{\ast}},$ respectively, whereas for
$x\in X$ and $x^{\ast}\in X^{\ast}$ we denote $\mathcal{N}_{x}:=x+\mathcal{N}%
_{X}$ and $\mathcal{N}_{x^{\ast}}:=x^{\ast}+\mathcal{N}_{X^{\ast}}.$

Given a nonempty set $A$ in $X$ (or in $X^{\ast}$), by $\operatorname*{co}A$,
$\operatorname*{cone}A$, and $\operatorname*{aff}A$ we denote the \emph{convex
hull,} the \emph{conic hull}, and\emph{ }the \emph{affine hull }of $A$,
respectively. Moreover, $\operatorname*{int}A$ is the \emph{interior} of $A$,
and $\operatorname*{cl}A$ and $\overline{A}$ are indistinctly used for
denoting the \emph{closure }of $A$ (weak*\emph{-}$\emph{closure}$ if $A\subset
X^{\ast}$). Thus, $\overline{\operatorname*{co}}A:=\operatorname*{cl}%
(\operatorname*{co}A)$, $\overline{\operatorname*{aff}}A:=\operatorname*{cl}%
(\operatorname*{aff}A),$ etc. We use $\operatorname*{ri}A$ to denote the
(topological) \emph{relative interior }of $A$ (i.e., the interior of $A$ in
the topology relative to $\overline{\operatorname*{aff}}A$). Associated with a
nonempty subset $A$ of $X,$ we consider the (one-sided)\emph{\ polar }and
the\emph{ orthogonal }of\emph{ }$A$ defined, respectively, by
\[
A^{\circ}:=\left\{  x^{\ast}\in X^{\ast}\mid\langle x^{\ast},x\rangle
\leq1\text{ for all }x\in A\right\}  ,
\]
and%
\[
A^{\perp}:=\left\{  x^{\ast}\in X^{\ast}\mid\langle x^{\ast},x\rangle=0\text{
for all }x\in A\right\}  .
\]
For $B\subset X$ and $\Omega\subset\mathbb{R}$ we define
\begin{align}
A+B  &  :=\{a+b\mid a\in A,\text{ }b\in B\},\nonumber\\
\Omega A  &  :=\{\lambda a\mid\lambda\in\Omega,\ a\in A\}, \label{Producto}%
\end{align}
with the conventions
\begin{equation}
\Omega\emptyset=A+\emptyset=\emptyset. \label{ConvProd}%
\end{equation}
Throughout\ the paper we apply the following two lemmas concerning product
(\ref{Producto}). Let us introduce the set%
\[
\Delta_{k}:=\{(\lambda_{1},\cdots,\lambda_{k})\mid\lambda_{i}>0,\sum
\nolimits_{i=1}^{k}\lambda_{i}=1\}.
\]

\begin{lem}
\label{Producto1}If $\Omega$ is a compact interval such that $0\notin\Omega$,
then%
\[
\Omega(\overline{\operatorname*{co}}A)=\overline{\operatorname*{co}}(\Omega
A).
\]

\end{lem}

\begin{dem}
The conclusion follows from the algebraic equality $\Omega(\operatorname*{co}%
A)=\operatorname*{co}(\Omega A)\ $together with its consequence
\[
\Omega(\overline{\operatorname*{co}}A)=\operatorname*{cl}(\Omega
(\operatorname*{co}A))=\overline{\operatorname*{co}}(\Omega A),
\]
both being true thanks to the condition $0\notin\Omega.$
\end{dem}

\begin{lem}
\label{Producto2}Suppose that $(\Lambda_{\varepsilon})_{\varepsilon>0}$ is a
non-increasing family of closed sets in $\mathbb{R}$ (i.e. $\varepsilon
^{\prime}<\varepsilon$ implies $\Lambda_{\varepsilon^{\prime}}\subset
\Lambda_{\varepsilon})$ such that $\bigcap\nolimits_{\varepsilon>0}%
\Lambda_{\varepsilon}=\{1\}.$ Let $(A_{\varepsilon})_{\varepsilon>0}$ be
another non-increasing family of closed sets in $X$ (or in $X^{\ast}).$ Then%
\[
\bigcap\nolimits_{\varepsilon>0}\Lambda_{\varepsilon}A_{\varepsilon}%
=\bigcap\nolimits_{\varepsilon>0}A_{\varepsilon}.
\]

\end{lem}

\begin{dem}
The inclusion "$\supset$" is obvious as $1\in\Lambda_{\varepsilon}$ entails
$A_{\varepsilon}\subset\Lambda_{\varepsilon}A_{\varepsilon}.$ Let us prove the
inclusion \textquotedblleft$\subset$\textquotedblright.

If $a\in\bigcap\nolimits_{\varepsilon>0}\Lambda_{\varepsilon}A_{\varepsilon}$,
we have $a=\lambda_{\varepsilon}a_{\varepsilon}$ with $\lambda_{\varepsilon
}\in\Lambda_{\varepsilon}$ and $a_{\varepsilon}\in A_{\varepsilon},$ for all
$\varepsilon>0.$ Since $\lim_{\varepsilon\downarrow0}\lambda_{\varepsilon}=1$
we can assume that $\lambda_{\varepsilon}>0.$ For each $\varepsilon>0$ the net
$(a_{\delta})_{0<\delta<\varepsilon}$ is contained in $A_{\varepsilon},$ and
due to the closedness of $A_{\varepsilon},$ $\lim_{\delta\downarrow0}%
a_{\delta}=\lim_{\delta\downarrow0}(\lambda_{\delta})^{-1}a=a\in
A_{\varepsilon},$ and we are done.
\end{dem}

We say that a function $\varphi:X\longrightarrow\overline{\mathbb{R}}$ is
proper if its \emph{(effective)} \emph{domain,} $\operatorname*{dom}%
\varphi:=\{x\in X\mid\varphi(x)<+\infty\}$, is nonempty. We say that $\varphi$
is \emph{convex }(\emph{lower semicontinuous }or\emph{ lsc}, for short,
respectively) if its \emph{epigraph,} $\operatorname*{epi}\varphi
:=\{(x,\lambda)\in X\times\mathbb{R}\mid\varphi(x)\leq\lambda\},$ is convex
(closed, respectively). We shall denote by $\Lambda(X)$ the set of all the
proper convex functions defined on $X$ and by $\Gamma_{0}(X)$ the functions in
$\Lambda(X)$ which are lsc. The \emph{lsc envelope }of $\varphi$ is the
function $\operatorname*{cl}\varphi$ such that $\operatorname*{epi}%
(\operatorname*{cl}\varphi)=\operatorname*{cl}(\operatorname*{epi}\varphi
)$.\emph{ }

If $\varepsilon\geq0$, the $\varepsilon$-\emph{subdifferential} of $\varphi$
at a point $x$ where $\varphi(x)$ is finite is the weak*-closed convex set
\[
\partial_{\varepsilon}\varphi(x):=\{x^{\ast}\in X^{\ast}\mid\varphi
(y)-\varphi(x)\geq\langle x^{\ast},y-x\rangle-\varepsilon\text{ for all }y\in
X\}.
\]
If $\varphi(x)\notin\mathbb{R}$, then we set $\partial_{\varepsilon}%
\varphi(x):=\emptyset$. In particular, for $\varepsilon=0$ we get the
\emph{Fenchel subdifferential} of $\varphi$ at $x,$ $\partial\varphi
(x):=\partial_{0}\varphi(x)$. When $\partial\varphi(x)\neq\emptyset,$ we know
that
\begin{equation}
\varphi(x)=(\operatorname*{cl}\varphi)(x)\text{ and }\partial\varphi
(x)=\partial(\operatorname*{cl}\varphi)(x)\text{.} \label{noemptys}%
\end{equation}

The \emph{support function,} the \emph{indicator function} and the
\emph{Minkowski gauge }of $A\subset X$ are, respectively, defined as
\begin{equation}
\sigma_{A}(x^{\ast}):=\sup\{\langle x^{\ast},a\rangle\mid a\in A\},\text{ for
}x^{\ast}\in X^{\ast},\text{ and }\sigma_{\emptyset}\equiv-\infty,
\label{support}%
\end{equation}%
\[
\mathrm{I}_{A}(x):=0,\text{ if }x\in A;\text{ }+\infty,\text{ if }x\in
X\setminus A,
\]
and
\[
p_{A}(x):=\inf\{\lambda\geq0:x\in\lambda A\},
\]
with $\inf\emptyset=+\infty.$ If $U\in\mathcal{N}_{X}$ (or $U\in
\mathcal{N}_{X^{\ast}}$), then $p_{U}$ is a continuous seminorm and
\begin{equation}
U=\{x\in X\mid p_{U}(x)\leq1\}. \label{gauge}%
\end{equation}
(Remember that $U$ is closed.)

If $A$ is convex, $x\in X$ and $\varepsilon\geq0$, we define the $\varepsilon
$\emph{-normal set }to $A$\ at $x$ as
\[
\mathrm{N}_{A}^{\varepsilon}(x):=\{x^{\ast}\in X^{\ast}\mid\langle x^{\ast
},y-x\rangle\leq\varepsilon\text{ for all }y\in A\},\text{ if \ }x\in A;\text{
\ }\emptyset,\text{ otherwise.}%
\]
If $\varepsilon=0$, we omit the reference to $\varepsilon$ and write
$\mathrm{N}_{A}(x).$

We shall frequently use the following well-known relation, for any set
$A\subset X$%
\begin{equation}
\bigcap\limits_{U\in\mathcal{N}_{X}}\operatorname*{cl}\left(  A+U\right)
=\operatorname*{cl}\left(  A\right)  . \label{ar}%
\end{equation}
The same relation is true for $A\subset X^{\ast}$ replacing $\mathcal{N}_{X}$
by $\mathcal{N}_{X^{\ast}}.$

\begin{lem}
\label{tec1}Let $A$ be a set in $X^{\ast}.$\ Then for every $x\in X$%
\[
\bigcap\limits_{L\in\mathcal{F(}x\mathcal{)}}\operatorname*{cl}\left(
A+L^{\perp}\right)  =\operatorname*{cl}\left(  A\right)  .
\]
If $X$ is a normed space, then
\[
\bigcap\limits_{\varepsilon>0}\operatorname*{cl}\left(  A+\varepsilon B_{\ast
}\right)  =\operatorname*{cl}(A),
\]
where $B_{\ast}$ is the closed unit ball in $X^{\ast}.$
\end{lem}

\begin{dem}
Due to the fact that the family of sets
\[
\{x^{\ast}\in X^{\ast}\mid\left\vert \left\langle x_{i},x^{\ast}\right\rangle
\right\vert \leq\delta,\ i\in\overline{1,k}\},\ \ x_{i}\in X,\ i\in
\overline{1,k},\ \ \delta>0\}
\]
is a basis of $\theta$-neighborhoods\ of the weak*-topology, and observing
that\textbf{ }%
\[
\{x_{i}\in X,\ i\in\overline{1,k};\ x\}^{\perp}\subset\{x^{\ast}\in X^{\ast
}\mid\left\vert \left\langle x_{i},x^{\ast}\right\rangle \right\vert
\leq\delta,\ i\in\overline{1,k}\},
\]
(\ref{ar}) allows us to write
\[
\operatorname*{cl}\left(  A\right)  \subset\bigcap\limits_{L\in\mathcal{F(}%
x\mathcal{)}}\operatorname*{cl}\left(  A+L^{\perp}\right)  \subset
\bigcap\limits_{U\in\mathcal{N}}\operatorname*{cl}\left(  A+U\right)
=\operatorname*{cl}\left(  A\right)  .
\]

Now, assume that\ $X$ is a normed space. Since the topology on the norm in
$X^{\ast}$ is stronger than every compatible topology in\ $X^{\ast},$ we get,
also using (\ref{ar}) in $X^{\ast},$%
\[
\operatorname*{cl}\left(  A\right)  \subset\bigcap\limits_{\varepsilon
>0}\operatorname*{cl}\left(  A+\varepsilon B_{\ast}\right)  \subset
\bigcap\limits_{U\in\mathcal{N}_{X^{\ast}}}\operatorname*{cl}\left(
A+U\right)  =\operatorname*{cl}\left(  A\right)  .
\]

\end{dem}

\section{The case of Banach spaces\label{Sect2b}}

Our purpose in this section is to provide a first direct consequence of
formula\ (\ref{formula using epsilon-subdifferential}), which yields a
free-continuity extension of (\ref{Valadier-like formula}) in a Banach space
$(X,\left\Vert \cdot\right\Vert ).$ This is accomplished via a straightforward
application of the Brøndsted-Rockafellar theorem. We also illustrate by
examples the validity limitation of this approach outside the Banach setting.
We need the following previous simple lemma.

\begin{lem}
\label{100}Given a convex function\ $\varphi:X\rightarrow\overline{\mathbb{R}%
},$ for every $x\in\operatorname*{dom}\varphi$ we have
\[
\bigcap_{L\in\mathcal{F}(x)}\partial(\varphi+\mathrm{I}_{L})(x)=\partial
\varphi(x).
\]

\end{lem}

For a function $\varphi\in\Gamma_{0}(X),$ we recall here the mapping
$\breve{\partial}^{\varepsilon}\!\varphi:X\rightrightarrows X^{\ast}$ defined
for $x\in\operatorname*{dom}\varphi$ as (see, e.g., \cite[(7) in p.
1268]{Lassonde})%
\begin{equation}
\breve{\partial}^{\varepsilon}\!\varphi(x):=\left\{  y^{\ast}\in X^{\ast
}\left\vert
\begin{tabular}
[c]{l}%
$\exists y\in B_{\varepsilon}(x),\text{ such that }\left\vert \varphi
(y)-\varphi(x)\right\vert \leq\varepsilon,$\\
$y^{\ast}\in\partial\varphi(y)\text{ and }\left\vert \langle y^{\ast
},y-x\rangle\right\vert \leq\varepsilon\}$%
\end{tabular}
\ \ \ \ \ \right.  \right\}  , \label{enlargement}%
\end{equation}
which constitutes an enlargement of the Fenchel subdifferential, and that can
also be written as
\[
\breve{\partial}^{\varepsilon}\!\varphi(x)=\bigcup_{y\in B_{\varphi
}(x,\varepsilon)}\partial\varphi(y)\cap\mathrm{N}_{\{y\}}^{\varepsilon}(x),
\]
where
\[
B_{\varphi}(x,\varepsilon):=\{y\in B_{\varepsilon}(x)\mid\left\vert
\varphi(y)-\varphi(x)\right\vert \leq\varepsilon\}.
\]
When $x\notin\operatorname*{dom}\varphi$ we adopt the natural convention
$\breve{\partial}^{\varepsilon}\!\varphi(x):=\emptyset.$ It is clear that
$\partial\varphi(x)\subset\breve{\partial}^{\varepsilon}\!\varphi(x)$ and,
moreover, it can be easily checked that $\cap_{\varepsilon>0}\breve{\partial
}^{\varepsilon}\!\varphi(x)=\partial\varphi(x).$ Let us emphasize here that
the elements of $\breve{\partial}^{\varepsilon}\!\varphi(x)$ are exact Fenchel
subgradients at nearby points to $x$, not $\varepsilon$-subgradients like in
(\ref{formula using epsilon-subdifferential}) and (\ref{marco5}).

The following version of the Brøndsted-Rockafellar theorem (see, e.g.
\cite[Theorem 3.1.1]{ZalinescuBook}) is a key tool in our approach.

\begin{lem}
\label{Brondsted-Rockafellar theorem}Let $\varphi\in\Gamma_{0}(X)$ and
$x\in\operatorname*{dom}\varphi.$ Then, for every $\varepsilon>0$ and
$x^{\ast}\in\partial_{\varepsilon}\varphi(x),$ there exist $x_{\varepsilon}\in
X,$ $y_{\varepsilon}^{\ast}\in B_{\ast},$ and $\lambda_{\varepsilon}\in\left[
-1,1\right]  $ such that
\[
\left\Vert x_{\varepsilon}-x\right\Vert \leq\sqrt{\varepsilon},
\]%
\[
x_{\varepsilon}^{\ast}:=x^{\ast}+\sqrt{\varepsilon}(y_{\varepsilon}^{\ast
}+\lambda_{\varepsilon}x^{\ast})\in\partial\varphi(x_{\varepsilon}),
\]%
\[
\left\vert \left\langle x_{\varepsilon}^{\ast},x_{\varepsilon}-x\right\rangle
\right\vert \leq\varepsilon+\sqrt{\varepsilon},\text{ }\left\vert
\varphi(x_{\varepsilon})-\varphi(x)\right\vert \leq\varepsilon+\sqrt
{\varepsilon},
\]
which implies that $x_{\varepsilon}^{\ast}\in\partial_{2\varepsilon}%
\varphi(x)$ and $x_{\varepsilon}^{\ast}\in\breve{\partial}^{\varepsilon
+\sqrt{\varepsilon}}\!\varphi(x).$
\end{lem}

Now we give\ the main result in this section.

\begin{theo}
\label{Thm valadier in Banach}Let\ $f_{t}\in\Gamma_{0}(X),$ $t\in T,$ and
$f=\sup_{t\in T}f_{t}.$ Then for every $x\in X$
\begin{equation}
\partial f(x)=\bigcap_{\substack{\varepsilon>0\\L\in\mathcal{F}(x)}%
}\overline{\operatorname*{co}}\left\{  \bigcup\limits_{t\in T_{\varepsilon
}(x)}\breve{\partial}^{\varepsilon}\!f_{t}(x)+\mathrm{N}_{L\cap
\operatorname*{dom}f}(x)\right\}  , \label{FVB}%
\end{equation}
where $T_{\varepsilon}(x)$ and $\mathcal{F}(x)$ have been defined in
\emph{(\ref{T-epsilon})} and \emph{(\ref{Fdex}),} respectively.
\end{theo}

\begin{dem}
If $x\notin\operatorname*{dom}f$ the formula holds trivially since both sets
$\mathrm{N}_{L\cap\operatorname*{dom}f}(x)$ and $\partial f(x)$ are empty and
we apply (\ref{ConvProd}); hence, we suppose that $x\in\operatorname*{dom}f$.

\emph{ Step 1. }To prove\ the inclusion \textquotedblleft$\supset
$\textquotedblright\ we first show that, for $\varepsilon>0,$ we have
\begin{equation}
\bigcup\limits_{t\in T_{\varepsilon}(x)}\breve{\partial}^{\varepsilon}%
\!f_{t}(x)\subset\partial_{3\varepsilon}f(x). \label{101}%
\end{equation}
Indeed, fix $t\in T_{\varepsilon}(x)$ and pick\ a $y^{\ast}\in\breve{\partial
}^{\varepsilon}\!f_{t}(x).$ By (\ref{enlargement}), there exists $y\in
B_{\varepsilon}(x)$ with $\left\vert f_{t}(y)-f_{t}(x)\right\vert
\leq\varepsilon\ $such that $y^{\ast}\in\partial f_{t}(y)$ and $\left\vert
\langle y^{\ast},y-x\rangle\right\vert \leq\varepsilon.\ $Then for every $z\in
X$
\begin{align*}
f(z)-f(x)  &  \geq f_{t}(z)-f_{t}(x)-\varepsilon\geq f_{t}(z)-f_{t}%
(y)-2\varepsilon\\
&  \geq\langle y^{\ast},z-y\rangle-2\varepsilon=\langle y^{\ast}%
,z-x\rangle+\langle y^{\ast},x-y\rangle-2\varepsilon\\
&  \geq\langle y^{\ast},z-x\rangle-3\varepsilon,
\end{align*}
and we get $y^{\ast}\in\partial_{3\varepsilon}f(x).$

\emph{Step 2. }Now (\ref{enlargement}) implies, for every $L\in\mathcal{F}%
(x),$
\begin{align*}
\bigcup\limits_{t\in T_{\varepsilon}(x)}\breve{\partial}^{\varepsilon}%
\!f_{t}(x)+\mathrm{N}_{L\cap\operatorname*{dom}f}(x)  &  \subset
\partial_{3\varepsilon}f(x)+\mathrm{N}_{L\cap\operatorname*{dom}f}(x)\\
&  \subset\partial_{3\varepsilon}(f+\mathrm{I}_{L\cap\operatorname*{dom}%
f})(x)=\partial_{3\varepsilon}(f+\mathrm{I}_{L})(x),
\end{align*}
so that, by convexifying and intersecting over $\varepsilon$ and $L,$ and
applying Lemma \ref{100},
\begin{align*}
\bigcap_{\substack{\varepsilon>0\\L\in\mathcal{F}(x)}}\overline
{\operatorname*{co}}\left\{  \bigcup\limits_{t\in T_{\varepsilon}(x)}%
\breve{\partial}^{\varepsilon}\!f_{t}(x)+\mathrm{N}_{L\cap\operatorname*{dom}%
f}(x)\right\}   &  \subset\bigcap_{\substack{\varepsilon>0\\L\in
\mathcal{F}(x)}}\partial_{3\varepsilon}(f+\mathrm{I}_{L})(x)\\
&  =\bigcap_{L\in\mathcal{F}(x)}\partial(f+\mathrm{I}_{L})(x)\\
&  =\partial f(x).
\end{align*}

\emph{Step 3. }To prove the inclusion \textquotedblleft$\subset$%
\textquotedblright\ we pick an $x^{\ast}\in\partial_{\varepsilon}f_{t}(x)$ for
$\varepsilon>0$ and $t\in T_{\varepsilon}(x).$ Then by Lemma
\ref{Brondsted-Rockafellar theorem} there are $x_{\varepsilon}\in X,$
$y_{\varepsilon}^{\ast}\in B_{\ast},$ and $\lambda_{\varepsilon}\in\left[
-1,1\right]  $ such that
\[
\left\Vert x_{\varepsilon}-x\right\Vert \leq\sqrt{\varepsilon},
\]%
\[
x_{\varepsilon}^{\ast}:=x^{\ast}+\sqrt{\varepsilon}(y_{\varepsilon}^{\ast
}+\lambda_{\varepsilon}x^{\ast})\in\partial f_{t}(x_{\varepsilon}),
\]%
\[
\left\vert \left\langle x_{\varepsilon}^{\ast},x_{\varepsilon}-x\right\rangle
\right\vert \leq\varepsilon+\sqrt{\varepsilon},\text{ }\left\vert
f_{t}(x_{\varepsilon})-f_{t}(x)\right\vert \leq\varepsilon+\sqrt{\varepsilon
},
\]
and
\[
x_{\varepsilon}^{\ast}\in\breve{\partial}^{\varepsilon+\sqrt{\varepsilon}%
}\!f_{t}(x).
\]
Assuming without loss of generality that $1\leq2(1-\sqrt{\varepsilon})$ (i.e.,
$\varepsilon<1/4$)$,$ we have%
\[
x^{\ast}\in\frac{1}{1+\lambda_{\varepsilon}\sqrt{\varepsilon}}\breve{\partial
}^{\varepsilon+\sqrt{\varepsilon}}f_{t}(x)+\frac{\sqrt{\varepsilon}}%
{1+\lambda_{\varepsilon}\sqrt{\varepsilon}}B_{\ast}\subset\Lambda
_{\varepsilon}\breve{\partial}^{\varepsilon+\sqrt{\varepsilon}}f_{t}%
(x)+2\sqrt{\varepsilon}B_{\ast},
\]
where $\Lambda_{\varepsilon}:=\left[  \frac{1}{1+\sqrt{\varepsilon}},\frac
{1}{1-\sqrt{\varepsilon}}\right]  .$ Hence, taking into account Lemma
\ref{tec1}, formula (\ref{formula using epsilon-subdifferential}) leads us to
\begin{align}
\partial f(x)  &  \subset\bigcap_{\substack{\varepsilon>0\\L\in\mathcal{F}%
(x)}}\overline{\operatorname*{co}}\left\{  \bigcup\limits_{t\in T_{\varepsilon
}(x)}\Lambda_{\varepsilon}\breve{\partial}^{\varepsilon+\sqrt{\varepsilon}%
}\!f_{t}(x)+2\sqrt{\varepsilon}B_{\ast}+\mathrm{N}_{L\cap\operatorname*{dom}%
f}(x)\right\} \nonumber\\
&  \subset\bigcap_{\delta>0}\bigcap_{\substack{\varepsilon>0\\L\in
\mathcal{F}(x)}}\overline{\operatorname*{co}}\left\{  \bigcup\limits_{t\in
T_{\varepsilon}(x)}\Lambda_{\varepsilon}\breve{\partial}^{\varepsilon
+\sqrt{\varepsilon}}\!f_{t}(x)+\delta B_{\ast}+\mathrm{N}_{L\cap
\operatorname*{dom}f}(x)\right\} \nonumber\\
&  =\bigcap_{\substack{\varepsilon>0\\L\in\mathcal{F}(x)}}\bigcap_{\delta
>0}\overline{\operatorname*{co}}\left\{  \bigcup\limits_{t\in T_{\varepsilon
}(x)}\Lambda_{\varepsilon}\breve{\partial}^{\varepsilon+\sqrt{\varepsilon}%
}\!f_{t}(x)+\delta B_{\ast}+\mathrm{N}_{L\cap\operatorname*{dom}f}(x)\right\}
\nonumber\\
&  =\bigcap_{\substack{\varepsilon>0\\L\in\mathcal{F}(x)}}\overline
{\operatorname*{co}}\left\{  \bigcup\limits_{t\in T_{\varepsilon}(x)}%
\Lambda_{\varepsilon}\breve{\partial}^{\varepsilon+\sqrt{\varepsilon}}%
\!f_{t}(x)+\mathrm{N}_{L\cap\operatorname*{dom}f}(x)\right\}  .
\label{marco57}%
\end{align}

\emph{Step 4. }Taking into account now that $\mathrm{N}_{L\cap
\operatorname*{dom}f}(x)$ is a cone we have
\[
\Lambda_{\varepsilon}\breve{\partial}^{\varepsilon+\sqrt{\varepsilon}}%
\!f_{t}(x)+\mathrm{N}_{L\cap\operatorname*{dom}f}(x)=\Lambda_{\varepsilon
}\left(  \breve{\partial}^{\varepsilon+\sqrt{\varepsilon}}\!f_{t}%
(x)+\mathrm{N}_{L\cap\operatorname*{dom}f}(x)\right)  ,
\]
and (\ref{marco57}) yields%
\begin{align}
\partial f(x)  &  \subset\bigcap_{\substack{\varepsilon>0\\L\in\mathcal{F}%
(x)}}\overline{\operatorname*{co}}\left\{  \bigcup\limits_{t\in T_{\varepsilon
}(x)}\Lambda_{\varepsilon}\left(  \breve{\partial}^{\varepsilon+\sqrt
{\varepsilon}}\!f_{t}(x)+\mathrm{N}_{L\cap\operatorname*{dom}f}(x)\right)
\right\} \nonumber\\
&  =\bigcap_{\substack{\varepsilon>0\\L\in\mathcal{F}(x)}}\overline
{\operatorname*{co}}\left\{  \Lambda_{\varepsilon}\left(  \bigcup\limits_{t\in
T_{\varepsilon}(x)}\breve{\partial}^{\varepsilon+\sqrt{\varepsilon}}%
\!f_{t}(x)+\mathrm{N}_{L\cap\operatorname*{dom}f}(x)\right)  \right\}
\label{marco61}\\
&  =\bigcap_{\substack{\varepsilon>0\\L\in\mathcal{F}(x)}}\Lambda
_{\varepsilon}\overline{\operatorname*{co}}\left\{  \bigcup\limits_{t\in
T_{\varepsilon}(x)}\breve{\partial}^{\varepsilon+\sqrt{\varepsilon}}%
\!f_{t}(x)+\mathrm{N}_{L\cap\operatorname*{dom}f}(x)\right\} \nonumber\\
&  =\bigcap_{\substack{\varepsilon>0\\L\in\mathcal{F}(x)}}\overline
{\operatorname*{co}}\left\{  \bigcup\limits_{t\in T_{\varepsilon}(x)}%
\breve{\partial}^{\varepsilon+\sqrt{\varepsilon}}\!f_{t}(x)+\mathrm{N}%
_{L\cap\operatorname*{dom}f}(x)\right\} \nonumber
\end{align}
where the second and the third equalities come from Lemma \ref{Producto1} and
Lemma \ref{Producto2}, respectively. Hence,%
\begin{align}
\partial f(x)  &  \subset\bigcap_{\substack{\varepsilon>0\\L\in\mathcal{F}%
(x)}}\overline{\operatorname*{co}}\left\{  \bigcup\limits_{t\in T_{\varepsilon
+\sqrt{\varepsilon}}(x)}\breve{\partial}^{\varepsilon+\sqrt{\varepsilon}%
}\!f_{t}(x)+\mathrm{N}_{L\cap\operatorname*{dom}f}(x)\right\} \nonumber\\
&  =\bigcap_{\substack{\eta>0\\L\in\mathcal{F}(x)}}\overline
{\operatorname*{co}}\left\{  \bigcup\limits_{t\in T_{\varepsilon}(x)}%
\breve{\partial}^{\eta}\!f_{t}(x)+\mathrm{N}_{L\cap\operatorname*{dom}%
f}(x)\right\}  .\nonumber
\end{align}

\end{dem}

\begin{rem}
\label{remb}\emph{Observe that formula (\ref{FVB})\ remains valid if instead
of }$\breve{\partial}^{\varepsilon}\!f_{t}(x)$\emph{ we use the larger set }%
\begin{equation}
\overset{\smallfrown}{\partial^{\varepsilon}}f_{t}(x):=\left\{  y^{\ast}\in
X^{\ast}\left\vert
\begin{tabular}
[c]{l}%
$\exists y\in X\text{, \emph{such that} }\left\vert f_{t}(y)-f_{t}%
(x)\right\vert \leq\varepsilon,$\\
$y^{\ast}\in\partial f_{t}(y)\text{ \emph{and} }\left\vert \langle y^{\ast
},y-x\rangle\right\vert \leq\varepsilon\}$%
\end{tabular}
\ \ \right.  \right\}  . \label{be}%
\end{equation}
\emph{The reason is that this set is also included in }$\partial
_{3\varepsilon}f(x)$\emph{ (see the proof of (\ref{101})).}
\end{rem}

If $X$ is finite-dimensional, $\mathrm{N}_{\operatorname*{dom}f}%
(x)\subset\mathrm{N}_{L\cap\operatorname*{dom}f}(x)$ for all $L\in
\mathcal{F}(x)$ and (\ref{FVB}) collapses to
\begin{equation}
\partial f(x)=\bigcap_{\varepsilon>0}\overline{\operatorname*{co}}\left\{
\bigcup\limits_{t\in T_{\varepsilon}(x)}\breve{\partial}^{\varepsilon}%
\!f_{t}(x)+\mathrm{N}_{\operatorname*{dom}f}(x)\right\}  . \label{FVB1}%
\end{equation}

In the following corollary we use an alternative enlargement of the Fenchel
subdifferential, which involves both subgradients at nearby points $y$ and
$\varepsilon$-subgradients at the nominal point $x.$ This provides a
characterization of $\partial f(x)$ which is at the same time a refinement of
formulas (\ref{formula using epsilon-subdifferential}) and (\ref{FVB}). For
instance, the $\varepsilon$-subgradients of the $f_{t}$'s involved in
(\ref{formula using epsilon-subdifferential}) can be chosen in the range of
$\partial f_{t}$ at sufficiently $f_{t}$-graphically close points to $x.$

\begin{cor}
\label{cor1}Under\ the assumptions of Theorem
\emph{\ref{Thm valadier in Banach}} we have, for every $x\in X,$
\[
\partial f(x)=\bigcap_{\substack{\varepsilon>0\\L\in\mathcal{F}(x)}%
}\overline{\operatorname*{co}}\left\{  \bigcup\limits_{t\in T_{\varepsilon
}(x)}\widehat{\partial}^{\varepsilon}f_{t}(x)+\mathrm{N}_{L\cap
\operatorname*{dom}f}(x)\right\}  ,
\]
where
\[
\widehat{\partial}^{\varepsilon}f_{t}(x):=\left\{  y^{\ast}\in X^{\ast
}\left\vert
\begin{tabular}
[c]{l}%
$\exists y\in B_{\varepsilon}(x)\text{ such that }\left\vert f_{t}%
(y)-f_{t}(x)\right\vert \leq\varepsilon$\\
and $y^{\ast}\in\partial f_{t}(y)\cap\partial_{2\varepsilon}f_{t}(x)$%
\end{tabular}
\ \right.  \right\}  .
\]

\end{cor}

\begin{dem}
Let us first check that
\[
\breve{\partial}^{\varepsilon}\!f_{t}(x)\subset\widehat{\partial}%
^{\varepsilon}f_{t}(x).
\]
If $x^{\ast}\in\breve{\partial}^{\varepsilon}\!f_{t}(x),$ then $x^{\ast}%
\in\partial f_{t}(y)$ for some $y\in B_{\varepsilon}(x)$ such that $\left\vert
f_{t}(y)-f_{t}(x)\right\vert \leq\varepsilon$ and $\left\vert \left\langle
x^{\ast},y-x\right\rangle \right\vert \leq\varepsilon.$ So, for every
$z\in\operatorname*{dom}f_{t},$%
\begin{align*}
\left\langle x^{\ast},z-x\right\rangle  &  =\left\langle x^{\ast
},z-y\right\rangle +\left\langle x^{\ast},y-x\right\rangle \\
&  \leq f_{t}(z)-f_{t}(y)+\varepsilon\\
&  \leq f_{t}(z)-f_{t}(x)+2\varepsilon,
\end{align*}
and we deduce that $x^{\ast}\in\partial_{2\varepsilon}f_{t}(x);$ that is,
$x^{\ast}\in\widehat{\partial}^{\varepsilon}f_{t}(x).$ Then, according to
Theorem \ref{Thm valadier in Banach} we get the direct inclusion
\textquotedblleft$\subset$\textquotedblright.

To check the converse inclusion \textquotedblleft$\supset$\textquotedblright%
\ we prove now that
\begin{equation}
\bigcup\limits_{t\in T_{\varepsilon}(x)}\widehat{\partial}^{\varepsilon}%
f_{t}(x)\subset\partial_{3\varepsilon}f(x). \label{marco201}%
\end{equation}
To this aim, take $x^{\ast}\in\bigcup\limits_{t\in T_{\varepsilon}%
(x)}\widehat{\partial}^{\varepsilon}f_{t}(x).$ Then, there is some $t_{0}\in
T_{\varepsilon}(x)$ and $y\in B_{\varepsilon}(x)$ such that $x^{\ast}%
\in\partial f_{t_{0}}(y)\cap\partial_{2\varepsilon}f_{t_{0}}(x),$ with
$\left\vert f_{t_{0}}(y)-f_{t_{0}}(x)\right\vert \leq\varepsilon.$ Then for
every $z\in\operatorname*{dom}f$
\begin{align*}
\left\langle x^{\ast},z-x\right\rangle  &  \leq\left\langle x^{\ast
},z-y\right\rangle +\left\langle x^{\ast},y-x\right\rangle \\
&  \leq f_{t_{0}}(z)-f_{t_{0}}(y)+f_{t_{0}}(y)-f_{t_{0}}(x)+2\varepsilon\\
&  =f_{t_{0}}(z)-f_{t_{0}}(x)+2\varepsilon\\
&  \leq f(z)-f(x)+3\varepsilon,
\end{align*}
that is, $x^{\ast}\in\partial_{3\varepsilon}f(x).$

Thanks to (\ref{marco201}) we write
\begin{align*}
\bigcap_{\substack{\varepsilon>0\\L\in\mathcal{F}(x)}}\overline
{\operatorname*{co}}\left\{  \bigcup\limits_{t\in T_{\varepsilon}%
(x)}\widehat{\partial}^{\varepsilon}f_{t}(x)+\mathrm{N}_{L\cap
\operatorname*{dom}f}(x)\right\}   &  \subset\bigcap_{\substack{\varepsilon
>0\\L\in\mathcal{F}(x)}}\overline{\operatorname*{co}}\left\{  \partial
_{3\varepsilon}f(x)+\mathrm{N}_{L\cap\operatorname*{dom}f}(x)\right\} \\
&  \subset\bigcap_{\substack{\varepsilon>0\\L\in\mathcal{F}(x)}}\partial
_{3\varepsilon}(f+\mathrm{I}_{L})(x)=\partial f(x).
\end{align*}

\end{dem}

The following example illustrates\ the need for\ the lower semicontinuity
assumption in Theorem \ref{Thm valadier in Banach}, as well as for\ the
condition\ $\left\vert f_{t}(y)-f_{t}(x)\right\vert \leq\varepsilon.$

\begin{exam}
\emph{Let }$X=\mathbb{R}$\emph{ and consider the functions }$f_{1}%
,f_{2}:X\rightarrow\mathbb{R}$\emph{ as }%
\[
f_{1}(x)=\left\{
\begin{tabular}
[c]{ll}%
$+\infty,\text{ }$ & $\text{\emph{if} }x<0,$\\
$1,\text{ }$ & $\text{\emph{if} }x=0,$\\
$x,\text{ }$ & $\text{\emph{if} }x>0,$%
\end{tabular}
\ \right.  \text{ and }f_{2}(x)=\left\{
\begin{tabular}
[c]{ll}%
$-x,\text{ }$ & $\text{\emph{if} }x<0,$\\
$1,\text{ }$ & $\text{\emph{if} }x=0,$\\
$+\infty,\text{ }$ & $\text{\emph{if} }x>0.$%
\end{tabular}
\ \right.
\]
\emph{So,} $f=\max\{f_{1},f_{2}\}=1+\mathrm{I}_{\{0\}}\ $\emph{so that
}$\partial f(0)=\mathbb{R}$\emph{,} \emph{but} \emph{for small} $\varepsilon
>0$%
\[
\breve{\partial}^{\varepsilon}f_{1}(0)=\partial f_{1}(0)=\emptyset,\text{
}\breve{\partial}^{\varepsilon}f_{2}(0)=\partial f_{2}(0)=\emptyset,
\]
\emph{and Theorem \ref{Thm valadier in Banach} fails.}
\end{exam}

\begin{exam}
\emph{Let }$X$\emph{ be any infinite-dimensional Banach space, let }$g$\emph{
be a non-continuous linear mapping,} \emph{and define the functions }%
$f_{t}:X\rightarrow\mathbb{R}$\emph{, }$t\in T:=\left]  0,+\infty\right[
,$\emph{ as}%
\[
f_{t}(x):=tg(x).
\]
\emph{So, }$f:=\sup_{t\in T}f_{t}=\mathrm{I}_{[g\leq0]}$\emph{ and we get
}$\partial f(\theta)=\mathrm{N}_{[g\leq0]}(\theta)\neq\emptyset.$ \emph{But
}$\partial f_{t}\equiv t\partial g\equiv\emptyset,$ \emph{so that}
$\breve{\partial}^{\varepsilon}f_{t}(\theta)\equiv\emptyset,$ \emph{for all}
$t\in T,$ \emph{and} \emph{Theorem \ref{Thm valadier in Banach} fails again.}
\end{exam}

The following example (e.g., \cite{HanLopZal2008}) shows that it is necessary
to consider\ exact subdifferentials at nearby points.

\begin{exam}
\emph{Let }$f_{1},f_{2}:\mathbb{R}\rightarrow\mathbb{R}\cup\{+\infty\}$\emph{
be defined by }%
\[
f_{1}(x)=\left\{
\begin{array}
[c]{ll}%
-\sqrt{x}, & \text{\emph{if} }x\geq0,\\
+\infty, & \text{\emph{if} }x<0,
\end{array}
\right.  \text{ \emph{and} }f_{2}(x)=f_{1}(-x).
\]
\emph{Then,} $f:=\max\{f_{1},f_{2}\}=\mathrm{I}_{\mathbb{\{}0\mathbb{\}}}$
\emph{so that }$\partial f(0)=\mathbb{R}.$ \emph{But,} $\partial
f_{1}(0)=\partial f_{2}(0)=\emptyset$ \emph{and so, for every} $\varepsilon
>0,$ \emph{we have } $T_{\varepsilon}(\theta)=\{1,2\}$ \emph{and}
\[
\bigcap_{\varepsilon>0}\overline{\operatorname*{co}}\left\{  (\partial
f_{1}(0)\cup\partial f_{2}(0))+\mathrm{N}_{\operatorname*{dom}f}(0)\right\}
=\emptyset.
\]
\emph{So, in order to recover the whole set }$\partial f(0)$\emph{, we need to
consider the subdifferentials of }$f_{1}$\emph{\ and }$f_{2}$ \emph{at nearby
points. Indeed, for }$\varepsilon<1$ \emph{we have }%
\begin{align*}
\breve{\partial}^{\varepsilon}\!f_{1}(0)  &  =\left\{  y^{\ast}\in X^{\ast
}\left\vert
\begin{tabular}
[c]{l}%
$\exists\ \left\vert y\right\vert \leq\varepsilon$ \emph{such that}
$\left\vert f_{1}(y)\right\vert \leq\varepsilon,$\\
$y^{\ast}\in\partial f_{1}(y)\text{ \emph{and} }\left\vert y^{\ast
},y\right\vert \leq\varepsilon\}$%
\end{tabular}
\ \ \ \right.  \right\} \\
&  =\left\{  -\frac{1}{2\sqrt{y}}\mid0<y\leq\varepsilon^{2}\right\}  =\left]
-\infty,-\frac{1}{2\varepsilon}\right]  ,
\end{align*}
\emph{and, similarly,}
\[
\breve{\partial}^{\varepsilon}\!f_{2}(0)=\left[  \frac{1}{2\varepsilon
},+\infty\right[  .
\]
\emph{Since} $T_{\varepsilon}(0)=\{1,2\},$
\[
\bigcup\limits_{t\in T_{\varepsilon}(0)}\breve{\partial}^{\varepsilon}%
\!f_{t}(0)+\mathrm{N}_{\operatorname*{dom}f}(0)=\left]  -\infty,-\frac
{1}{2\varepsilon}\right]  \cup\left[  \frac{1}{2\varepsilon},+\infty\right[
+\mathbb{R}=\mathbb{R}.
\]

\end{exam}

We close this section with the following example showing that, in general,
Theorem \ref{Thm valadier in Banach} is not valid outside the Banach spaces setting.

\begin{exam}
\label{ExampleKlee}\emph{Let }$X$\emph{ be a lcs space such that there
exists\ a proper lsc convex function }$g\in\Gamma_{0}(X)$ \emph{having an
empty subdifferential everywhere. This is the case of\ some Fréchet spaces
(\cite{Rockafellar65}), and also of some non-complete normed spaces (actually
certain subspaces of }$l^{2}(\mathbb{N})$\emph{ \cite{Borwein85supportless}).}

\emph{According to \cite{Rockafellar65}, we may suppose that }$\theta
\in\operatorname*{dom}g$ \emph{and }$g(\theta)=0.$\emph{ For }$t\in T:=\left]
0,+\infty\right[  ,$\emph{ we define the function }$f_{t}\in\Gamma_{0}%
(X)$\emph{ as}%
\[
f_{t}(x):=tg(x),
\]
\emph{so that }$f:=\sup_{t\in T}f_{t}=\mathrm{I}_{[g\leq0]}.$\emph{ Since
}$\partial f_{t}\equiv t\partial f\equiv\emptyset$\emph{, the right-hand set
in (\ref{FVB}) is empty, whereas }%
\[
\partial f(\theta)=\mathrm{N}_{[g\leq0]}(\theta)\neq\emptyset.
\]

\end{exam}

\section{The case of locally convex spaces\label{Sect3}}

We give in this section a general free-continuity formula characterizing the
subdifferential of the supremum function in the setting of a lcs $X.$ In what
follows, $\mathcal{P}$ is the family of continuous seminorms in $X.$

We adapt the mapping introduced in (\ref{enlargement}) to our current setting:
given a convex function $\varphi:X\rightarrow\mathbb{R\cup\{+\infty\}}$ and
$x\in\operatorname*{dom}\varphi,$\ for $p\in\mathcal{P}$\ we denote
\begin{equation}
\breve{\partial}_{p}^{\varepsilon}\varphi(x):=\left\{  y^{\ast}\in X^{\ast
}\left\vert
\begin{array}
[c]{l}%
\exists y\in X\text{ with }p(y-x)\leq\varepsilon,\text{ such that\ }\\
\left\vert \varphi(y)-\varphi(x)\right\vert \leq\varepsilon,\text{ }y^{\ast
}\in\partial\varphi(y)\text{ and }\left\vert \langle y^{\ast},y-x\rangle
\right\vert \leq\varepsilon
\end{array}
\right.  \right\}  . \label{17}%
\end{equation}
Also now, when $x\notin\operatorname*{dom}\varphi$ we set $\breve{\partial
}_{p}^{\varepsilon}\!\varphi(x)=\emptyset.$

\begin{rem}
\label{remnobar}\emph{It is worth observing that }$\breve{\partial}%
_{p}^{\varepsilon}\varphi(x)$\emph{ remains the same if we replace the
inequality }$\left\vert \langle y^{\ast},y-x\rangle\right\vert \leq
\varepsilon$\emph{ by }$\langle y^{\ast},y-x\rangle\leq\varepsilon,$\emph{
thanks to the fact that }$y^{\ast}\in\partial\varphi(y)$\emph{ and
}$\left\vert \varphi(y)-\varphi(x)\right\vert \leq\varepsilon.$
\end{rem}

We give the main theorem of this section.

\begin{theo}
\label{Thm lcs1}Let\ $f_{t}\in\Gamma_{0}(X),$ $t\in T,$ and $f=\sup_{t\in
T}f_{t}.$ Then for every $x\in X$
\begin{equation}
\partial f(x)=\bigcap_{_{\substack{\varepsilon>0,\text{ }p\in\mathcal{P}%
\\L\in\mathcal{F}(x)}}}\overline{\operatorname*{co}}\left\{  \bigcup
\limits_{t\in T_{\varepsilon}(x)}\breve{\partial}_{p}^{\varepsilon}%
(f_{t}+\mathrm{I}_{\overline{L\cap\operatorname*{dom}f}})(x)\right\}  .
\label{marco1}%
\end{equation}

\end{theo}

\begin{dem}
If $x\notin\operatorname*{dom}f,$ then $x\notin\operatorname*{dom}f_{t}$ for
$t\in T_{\varepsilon}(x);$ hence, $x\notin\operatorname*{dom}(f_{t}%
+\mathrm{I}_{\overline{L\cap\operatorname*{dom}f}})$ and the sets in both
sides of (\ref{marco1}) are empty.

Now we assume that $x\in\operatorname*{dom}f.$ First we verify the inclusion
\textquotedblleft$\supset$\textquotedblright\ in (\ref{marco1}). We take
$L\in\mathcal{F}(x)$ and $y^{\ast}\in\breve{\partial}_{p}^{\varepsilon}%
(f_{t}+\mathrm{I}_{\overline{L\cap\operatorname*{dom}f}})(x)$ for some $t\in
T_{\varepsilon}(x).$ Let $y\in\overline{L\cap\operatorname*{dom}f}$ be such
that $p(y-x)\leq\varepsilon,$
\[
\left\vert (f_{t}+\mathrm{I}_{\overline{L\cap\operatorname*{dom}f}}%
)(y)-(f_{t}+\mathrm{I}_{\overline{L\cap\operatorname*{dom}f}})(x)\right\vert
=\left\vert f_{t}(y)-f_{t}(x)\right\vert \leq\varepsilon,
\]
$y^{\ast}\in\partial(f_{t}+\mathrm{I}_{\overline{L\cap\operatorname*{dom}f}%
})(y)$ and $\left\vert \langle y^{\ast},y-x\rangle\right\vert \leq
\varepsilon.$ Then, for all $z\in L\cap\operatorname*{dom}f,$\
\[
\left\langle y^{\ast},z-x\right\rangle =\left\langle y^{\ast},z-y\right\rangle
+\left\langle y^{\ast},y-x\right\rangle \leq f_{t}(z)-f_{t}(y)+\varepsilon\leq
f(z)-f(x)+3\varepsilon,
\]
and we get $y^{\ast}\in\partial_{3\varepsilon}(f+\mathrm{I}_{L})(x);$ that is,
$\breve{\partial}_{p}^{\varepsilon}(f_{t}+\mathrm{I}_{\overline{L\cap
\operatorname*{dom}f}})(x)\subset\partial_{3\varepsilon}(f+\mathrm{I}%
_{L})(x).$ Thus, recalling Lemma \ref{100},%
\begin{align*}
\bigcap_{_{\substack{\varepsilon>0,\text{ }p\in\mathcal{P}\\L\in
\mathcal{F}(x)}}}\overline{\operatorname*{co}}\left\{  \bigcup\limits_{t\in
T_{\varepsilon}(x)}\breve{\partial}_{p}^{\varepsilon}(f_{t}+\mathrm{I}%
_{\overline{L\cap\operatorname*{dom}f}})(x)\right\}   &  \subset\bigcap
_{L\in\mathcal{F}(x)}\bigcap_{\varepsilon>0}\partial_{3\varepsilon
}(f+\mathrm{I}_{L})(x)\\
&  =\bigcap_{L\in\mathcal{F}(x)}\partial(f+\mathrm{I}_{L})(x)=\partial f(x).
\end{align*}

To prove the inclusion \textquotedblleft$\subset$\textquotedblright\ in
(\ref{marco1}) we shall proceed by steps.

\emph{Step 1. }Observe that, due to the following\ relations
\[
\partial_{\varepsilon}f_{t}(x)+\mathrm{N}_{L\cap\operatorname*{dom}%
f}(x)=\partial_{\varepsilon}f_{t}(x)+\mathrm{N}_{\overline{L\cap
\operatorname*{dom}f}}(x)\subset\partial_{\varepsilon}(f_{t}+\mathrm{I}%
_{\overline{L\cap\operatorname*{dom}f}})(x),
\]
formula (\ref{formula using epsilon-subdifferential}) implies that
\begin{equation}
\partial f(x)\subset\bigcap_{\substack{\varepsilon>0\\L\in\mathcal{F}%
(x)}}\overline{\operatorname*{co}}\left\{  \bigcup\limits_{t\in T_{\varepsilon
}(x)}\partial_{\varepsilon}(f_{t}+\mathrm{I}_{\overline{L\cap
\operatorname*{dom}f}})(x)\right\}  . \label{clc}%
\end{equation}

\emph{Step 2. }Now, we\ fix $p\in\mathcal{P}$, $U\in\mathcal{N}_{X^{\ast}},$
$L\in\mathcal{F}(x),$ and $\varepsilon>0.$ Let $\left\Vert \cdot\right\Vert $
be a norm in $L$ such that $p\leq\left\Vert \cdot\right\Vert $ in such a
subspace, and take into account that $L^{\ast}\simeq X^{\ast}/L^{\perp}$ so
that $\mathcal{N}_{X^{\ast}}+L^{\perp}$ is\ a basis of $\theta$-neighborhoods
in $X^{\ast}/L^{\perp}$ (see, e.g., \cite{Fabianetal}). Choose $\eta>0$ such
that
\begin{equation}
\eta+\sqrt{\eta}\leq\varepsilon\text{ and }\sqrt{\eta}B_{L^{\ast}}\subset
U+L^{\perp}, \label{a0}%
\end{equation}
where $B_{L^{\ast}}$ is the unit ball in $L^{\ast}$.

For the sake of notational brevity, we shall denote
\[
g_{t}:=f_{t}+\mathrm{I}_{\overline{L\cap\operatorname*{dom}f}},\text{ }t\in
T,
\]
where $L$ is our fixed subspace.

Take $x^{\ast}\in\partial_{\eta}g_{t}(x),$ with $t\in T_{\eta}(x).$ Then
$x_{\mid L}^{\ast}\in\partial_{\eta}g_{t_{\mid L}}(x),$ and by Proposition
\ref{Brondsted-Rockafellar theorem} we find $x_{\eta}\in L$ $(\subset X)$ and
$v_{\eta}^{\ast}\in\partial g_{t_{\mid L}}(x_{\eta})$ (hence, $x_{\eta}%
\in\overline{L\cap\operatorname*{dom}f})$, together with $\lambda_{\eta}%
\in\left[  -1,1\right]  $ and $y^{\ast}\in B_{L^{\ast}}$ (hence, $\sqrt{\eta
}y^{\ast}\in\sqrt{\eta}B_{L^{\ast}}\subset U+L^{\perp}$), such that
\begin{equation}
p(x_{\eta}-x)\leq\left\Vert x_{\eta}-x\right\Vert \leq\sqrt{\eta}%
\leq\varepsilon, \label{a1}%
\end{equation}%
\[
v_{\eta}^{\ast}:=x_{\mid L}^{\ast}+\sqrt{\eta}(y^{\ast}+\lambda_{\eta}x_{\mid
L}^{\ast})=(1+\lambda_{\eta}\sqrt{\eta})x_{\mid L}^{\ast}+\sqrt{\eta}y^{\ast
},
\]%
\[
\left\vert \left\langle v_{\eta}^{\ast},x_{\eta}-x\right\rangle \right\vert
\leq\eta+\sqrt{\eta}\leq\varepsilon,
\]
and%
\begin{equation}
\left\vert f_{t}(x_{\eta})-f_{t}(x)\right\vert =\left\vert g_{t}(x_{\eta
})-g_{t}(x)\right\vert \leq\eta+\sqrt{\eta}\leq\varepsilon. \label{a2}%
\end{equation}
\emph{ }

By the Hahn-Banach theorem we extend $v_{\eta}^{\ast}$ to an $x_{\eta}^{\ast
}\in X^{\ast};$ hence, in particular, we have $x_{\eta}^{\ast}\in\partial
g_{t}(x_{\eta})\ $and $\left\vert \langle x_{\eta}^{\ast},x_{\eta}%
-x\rangle\right\vert =\left\vert \left\langle v_{\eta}^{\ast},x_{\eta
}-x\right\rangle \right\vert \leq\varepsilon,$ which together with (\ref{a1})
and (\ref{a2}) ensures that
\begin{equation}
x_{\eta}^{\ast}\in\breve{\partial}_{p}^{\varepsilon}g_{t}(x). \label{am3}%
\end{equation}
Hence we have proved that%
\[
\partial_{\eta}g_{t}(x)\subset\breve{\partial}_{p}^{\varepsilon}%
g_{t}(x),\text{ for all }t\in T_{\eta}(x).
\]

\emph{Step 3}$.$ Moreover, for every $u\in L$ we obtain
\begin{align*}
\left\langle (1+\lambda_{\eta}\sqrt{\eta})x^{\ast},u-x\right\rangle  &
=\left\langle (1+\lambda_{\eta}\sqrt{\eta})x_{\mid L}^{\ast},u-x\right\rangle
\\
&  =\left\langle (1+\lambda_{\eta}\sqrt{\eta})x_{\mid L}^{\ast}-v_{\eta}%
^{\ast},u-x\right\rangle +\left\langle v_{\eta}^{\ast},u-x\right\rangle \\
&  =\left\langle -\sqrt{\eta}y^{\ast},u-x\right\rangle +\left\langle x_{\eta
}^{\ast},u-x\right\rangle \\
&  =\left\langle x_{\eta}^{\ast}-\sqrt{\eta}y^{\ast},u-x\right\rangle
\end{align*}
that is, taking into account that $\sqrt{\eta}y^{\ast}\in\sqrt{\eta}%
B_{L^{\ast}}\subset U+L^{\perp}$ (by (\ref{a0})), and using (\ref{am3}),%
\begin{align}
(1+\lambda_{\eta}\sqrt{\eta})\left\langle x^{\ast},u-x\right\rangle  &
=\left\langle x_{\eta}^{\ast}-\sqrt{\eta}y^{\ast},u-x\right\rangle \nonumber\\
&  \leq\sigma_{\breve{\partial}_{p}^{\varepsilon}g_{t}(x)+U+L^{\perp}}(u-x).
\label{19}%
\end{align}
Observe that (\ref{19}) can be extended to $u\in X$, i.e. we also have%
\begin{equation}
(1+\lambda_{\eta}\sqrt{\eta})\left\langle x^{\ast},u-x\right\rangle \leq
\sigma_{\breve{\partial}_{p}^{\varepsilon}g_{t}(x)+U+L^{\perp}}(u-x)\text{ for
all }u\in X, \label{19'}%
\end{equation}
due to the fact that, for any $u\in X$,
\[
\sigma_{\breve{\partial}_{p}^{\varepsilon}g_{t}(x)+U+L^{\perp}}(u-x)=\sup
\{\langle u^{\ast}+\ell^{\ast},u-x\rangle\mid u^{\ast}\in\breve{\partial}%
_{p}^{\varepsilon}g_{t}(x)+U,\ \ell^{\ast}\in L^{\perp}\},\text{ }%
\]
and the expression on the right-hand is $+\infty$ when $u-x\not \in L,$
equivalently $u\not \in L.$

\emph{Step 4. }Now observe that
\begin{equation}
\breve{\partial}_{p}^{\varepsilon}g_{t}(x)+L^{\perp}\subset\breve{\partial
}_{p}^{\varepsilon}g_{t}(x). \label{18}%
\end{equation}
Indeed, if $z^{\ast}\in\breve{\partial}_{p}^{\varepsilon}g_{t}(x)$ and
$u^{\ast}\in L^{\perp},$ then by (\ref{17}) there exists $z\in X$ such that
$p(z-x)\leq\varepsilon,$ $\left\vert g_{t}(z)-g_{t}(x)\right\vert
\leq\varepsilon,$ $z^{\ast}\in\partial g_{t}(z)$ and $\left\vert \langle
z^{\ast},z-x\rangle\right\vert \leq\varepsilon;$ in particular, $z\in
\operatorname*{dom}g_{t}\subset L.$ Then,
\begin{align*}
z^{\ast}+u^{\ast}  &  \in\partial g_{t}(z)+L^{\perp}\subset\partial
(g_{t}+\mathrm{I}_{L})(z)\\
&  =\partial(f_{t}+\mathrm{I}_{\overline{L\cap\operatorname*{dom}f}%
}+\mathrm{I}_{L})(z)=\partial g_{t}(z).
\end{align*}
As $\left\vert \langle z^{\ast}+u^{\ast},z-x\rangle\right\vert =\left\vert
\langle z^{\ast},z-x\rangle\right\vert \leq\varepsilon$ we deduce that
$z^{\ast}+u^{\ast}\in\breve{\partial}_{p}^{\varepsilon}g_{t}(x),$ yielding
(\ref{18}).

\emph{Step 5. }Consequently, using (\ref{19'}) and (\ref{18}), for all $z\in
X$ (as $t\in T_{\eta}(x)\subset T_{\varepsilon}(x))$%
\begin{align*}
(1+\lambda_{\eta}\sqrt{\eta})\left\langle x^{\ast},z-x\right\rangle  &
\leq\sigma_{\breve{\partial}_{p}^{\varepsilon}g_{t}(x)+U}(z-x)\\
&  \leq\sigma_{\bigcup\limits_{t\in T_{\eta}(x)}\breve{\partial}%
_{p}^{\varepsilon}g_{t}(x)+U}(z-x)\\
&  \leq\sigma_{\bigcup\limits_{t\in T_{\varepsilon}(x)}\breve{\partial}%
_{p}^{\varepsilon}g_{t}(x)+U}(z-x)\\
&  =\sigma_{\overline{\operatorname*{co}}\left\{  \bigcup\limits_{t\in
T_{\varepsilon}(x)}\breve{\partial}_{p}^{\varepsilon}g_{t}(x)+U\right\}
}(z-x);
\end{align*}
therefore
\[
(1+\lambda_{\eta}\sqrt{\eta})x^{\ast}\in\overline{\operatorname*{co}}\left\{
\bigcup\limits_{t\in T_{\varepsilon}(x)}\breve{\partial}_{p}^{\varepsilon
}g_{t}(x)+U\right\}
\]
and so, recalling that $x^{\ast}\in\partial_{\eta}g_{t}(x),$ with $t\in
T_{\eta}(x),$
\begin{equation}
\bigcup\limits_{t\in T_{\eta}(x)}\partial_{\eta}g_{t}(x)\subset\Lambda_{\eta
}\overline{\operatorname*{co}}\left\{  \bigcup\limits_{t\in T_{\varepsilon
}(x)}\breve{\partial}_{p}^{\varepsilon}g_{t}(x)+U\right\}  , \label{marco97}%
\end{equation}
where $\Lambda_{\eta}=\left[  \frac{1}{1+\sqrt{\eta}},\frac{1}{1-\sqrt{\eta}%
}\right]  .$

\emph{Step 6. }Now, from (\ref{marco97}) and Lemma \ref{Producto1} which
ensures that the set in the right-hand side above is closed and convex,%
\begin{align*}
\bigcap_{\substack{\delta>0\\M\in\mathcal{F}(x)}}\overline{\operatorname*{co}%
}\left\{  \bigcup\limits_{t\in T_{\delta}(x)}\partial_{\delta}(f_{t}%
+\mathrm{I}_{\overline{M\cap\operatorname*{dom}f}})(x)\right\}   &
\subset\overline{\operatorname*{co}}\left\{  \bigcup\limits_{t\in T_{\eta}%
(x)}\breve{\partial}_{p}^{\eta}g_{t}(x)\right\} \\
&  \subset\Lambda_{\eta}\overline{\operatorname*{co}}\left\{  \bigcup
\limits_{t\in T_{\varepsilon}(x)}\breve{\partial}_{p}^{\varepsilon}%
g_{t}(x)+U\right\}  ,
\end{align*}
so that, intersecting over $\varepsilon>0$ (recall (\ref{clc}) and Lemma
(\ref{Producto2})),%
\[
\partial f(x)\subset\bigcap_{\varepsilon>0}\Lambda_{\eta}\overline
{\operatorname*{co}}\left\{  \bigcup\limits_{t\in T_{\varepsilon}(x)}%
\breve{\partial}_{p}^{\varepsilon}g_{t}(x)+U\right\}  =\bigcap_{\varepsilon
>0}\overline{\operatorname*{co}}\left\{  \bigcup\limits_{t\in T_{\varepsilon
}(x)}\breve{\partial}_{p}^{\varepsilon}g_{t}(x)+U\right\}  ,
\]
which in turn gives us, by intersecting over $U$ (recall (\ref{ar})),%
\[
\partial f(x)\subset\bigcap_{\varepsilon>0}\overline{\operatorname*{co}%
}\left\{  \bigcup\limits_{t\in T_{\varepsilon}(x)}\breve{\partial}%
_{p}^{\varepsilon}g_{t}(x)\right\}  .
\]
Finally, the inclusion follows since $L$ and $p$ were arbitrarily chosen.
\end{dem}

\bigskip

\begin{rem}
\emph{(\ref{marco1}) gives an alternative formula of (\ref{FVB}) in the
setting of lcs spaces. If the underlying space is Banach, all that can be
derived from (\ref{marco1}), by using\ Proposition
\ref{Brondsted-Rockafellar theorem}, is }%
\[
\breve{\partial}_{p}^{\varepsilon}(f_{t}+\mathrm{I}_{\overline{L\cap
\operatorname*{dom}f}})(x)\subset\breve{\partial}_{p}^{\varepsilon}f_{t}%
(y_{1})+\mathrm{N}_{\overline{L\cap\operatorname*{dom}f}}(y_{2}),
\]
\emph{for some }$y_{1}$\emph{ and }$y_{2}$\emph{ close to }$x.$\emph{ In this
sense, we can not directly deduce (\ref{FVB}) from (\ref{marco1}) (in Banach
spaces).}
\end{rem}

\begin{rem}
\emph{Similar to the observation made in Remark \ref{remb}, we can
equivalently write formula (\ref{marco1}) as}
\[
\partial f(x)=\bigcap_{_{\substack{\varepsilon>0\\L\in\mathcal{F}(x)}%
}}\overline{\operatorname*{co}}\left\{  \bigcup\limits_{t\in T_{\varepsilon
}(x)}\overset{\smallfrown}{\partial^{\varepsilon}}\!(f_{t}+\mathrm{I}%
_{\overline{L\cap\operatorname*{dom}f}})(x)\right\}  ,
\]
\emph{where} $\overset{\smallfrown}{\partial^{\varepsilon}}$ \emph{is the
mapping introduced in (\ref{be}).}
\end{rem}

\begin{rem}
\emph{As in Corollary \ref{cor1}, we can show that }%
\[
\breve{\partial}_{p}^{\varepsilon}(f_{t}+\mathrm{I}_{\overline{L\cap
\operatorname*{dom}f}})(x)\subset\widehat{\partial}_{p}^{\varepsilon}%
(f_{t}+\mathrm{I}_{\overline{L\cap\operatorname*{dom}f}})(x),
\]
\emph{where for a function }$\varphi\in\Lambda(X)$ \emph{we denote\ }%
\begin{equation}
\widehat{\partial}_{p}^{\varepsilon}\varphi(x):=\left\{  y^{\ast}\in X^{\ast
}\left\vert
\begin{tabular}
[c]{l}%
$\exists y\in X\text{ \emph{such that} }p(y-x)\leq\varepsilon,$ $\left\vert
\varphi(y)-\varphi(x)\right\vert \leq\varepsilon$\\
\emph{and} $y^{\ast}\in\partial\varphi(y)\cap\partial_{2\varepsilon}%
\varphi(x)$%
\end{tabular}
\ \ \right.  \right\}  . \label{gorop}%
\end{equation}
\emph{As a consequence of (\ref{gorop}) we have that }%
\begin{equation}
\breve{\partial}_{p}^{\varepsilon}(f_{t}+\mathrm{I}_{\overline{L\cap
\operatorname*{dom}f}})(x)\subset\widehat{\partial}_{p}^{\varepsilon}%
(f_{t}+\mathrm{I}_{\overline{L\cap\operatorname*{dom}f}})(x)\subset
\partial_{2\varepsilon}(f_{t}+\mathrm{I}_{\overline{L\cap\operatorname*{dom}%
f}})(x). \label{cco}%
\end{equation}
\emph{Then similar arguments to those used in Corollary \ref{cor1} give rise
to }%
\[
\partial f(x)=\bigcap_{_{\substack{\varepsilon>0\\L\in\mathcal{F}(x)}%
}}\overline{\operatorname*{co}}\left\{  \bigcup\limits_{t\in T_{\varepsilon
}(x)}\widehat{\partial}_{p}^{\varepsilon}(f_{t}+\mathrm{I}_{\overline
{L\cap\operatorname*{dom}f}})(x)\right\}  .
\]

\end{rem}

\bigskip

Once again the intersection over the $L$'s in (\ref{marco1}) can be removed in
the finite-dimensional setting. In the following theorem we show that this is
also the case in a more general setting.

\begin{theo}
\label{sinl}Under the assumptions of Theorem \emph{\ref{Thm lcs1},}
if\ $\operatorname*{ri}(\operatorname*{dom}f)\neq\emptyset$ and $f_{\mid
\operatorname*{aff}(\operatorname*{dom}f)}$ is continuous on
$\operatorname*{ri}(\operatorname*{dom}f),$ then for all $x\in X$
\[
\partial f(x)=\bigcap_{_{\varepsilon>0,\text{ }p\in\mathcal{P}}}%
\overline{\operatorname*{co}}\left\{  \bigcup\limits_{t\in T_{\varepsilon}%
(x)}\breve{\partial}_{p}^{\varepsilon}(f_{t}+\mathrm{I}_{\overline
{\operatorname*{dom}f}})(x)\right\}  .
\]

\end{theo}

\begin{dem}
Fix $x\in\operatorname*{dom}f$ and choose $L\in\mathcal{F}(x)$ such that
$L\cap\operatorname*{ri}(\operatorname*{dom}f)\neq\emptyset;$ hence by the
accessibility lemma%
\begin{equation}
\overline{L\cap\operatorname*{dom}f}=L\cap\overline{\operatorname*{dom}f}.
\label{la}%
\end{equation}
Next, given $\varepsilon>0,$ $p\in\mathcal{P}$ and $t\in T_{\varepsilon}(x),$
we show that
\begin{equation}
\breve{\partial}_{p}^{\varepsilon}(f_{t}+\mathrm{I}_{\overline{L\cap
\operatorname*{dom}f}})(x)\subset\operatorname*{cl}\left(  \breve{\partial
}_{p}^{2\varepsilon}(f_{t}+\mathrm{I}_{\overline{\operatorname*{dom}f}%
})(x)+L^{\bot}\right)  . \label{li}%
\end{equation}
Take $y^{\ast}\in\breve{\partial}_{p}^{\varepsilon}(f_{t}+\mathrm{I}%
_{\overline{L\cap\operatorname*{dom}f}})(x)$ and let $y\in\overline
{L\cap\operatorname*{dom}f}$ ($\subset L$) be such that $p(y-x)\leq
\varepsilon,$ $\left\vert f_{t}(y)-f_{t}(x)\right\vert \leq\varepsilon,$
$\left\vert \langle y^{\ast},y-x\rangle\right\vert \leq\varepsilon,$ and
$y^{\ast}\in\partial(f_{t}+\mathrm{I}_{\overline{L\cap\operatorname*{dom}f}%
})(y).$

Denote $h_{t}:=f_{t}+\mathrm{I}_{\overline{\operatorname*{dom}f}},$ $t\in T.$
Since $\operatorname*{dom}f\subset\operatorname*{dom}f_{t}\cap\overline
{\operatorname*{dom}f}=\operatorname*{dom}h_{t}\subset\overline
{\operatorname*{dom}f},$ we have\ that
\[
\operatorname*{ri}(\operatorname*{dom}f)=\operatorname*{ri}%
(\operatorname*{dom}f_{t}\cap\overline{\operatorname*{dom}f}%
)=\operatorname*{ri}(\operatorname*{dom}h_{t})\text{ and }\operatorname*{aff}%
(\operatorname*{dom}h_{t})=\operatorname*{aff}(\operatorname*{dom}f).
\]
Consequently, because $h_{t\mid\operatorname*{aff}(\operatorname*{dom}g_{t}%
)}\leq f_{\mid\operatorname*{aff}(\operatorname*{dom}f)}$ the function
$h_{t\mid\operatorname*{aff}(\operatorname*{dom}g_{t})}$ is continuous on
$\operatorname*{ri}(\operatorname*{dom}h_{t})$ and so, according to
\cite[Theorem 15]{vv} (remember\ (\ref{la})),%
\[
y^{\ast}\in\partial(f_{t}+\mathrm{I}_{\overline{L\cap\operatorname*{dom}f}%
})(y)=\partial(h_{t}+\mathrm{I}_{L})(y)=\operatorname*{cl}(\partial
h_{t}(y)+L^{\bot}).
\]
Let $U_{0}\in\mathcal{N}_{X^{\ast}}$ be such that $\sigma_{U_{0}}%
(y-x)\leq\varepsilon.$ Then for every $U\in\mathcal{N}_{X^{\ast}}$ such that
$U\subset U_{0}$ there exists $z^{\ast}\in\partial h_{t}(y)$ such that
$y^{\ast}\in z^{\ast}+L^{\bot}+U$ and so, as $x,y\in L,$
\[
\left\vert \langle z^{\ast},y-x\rangle\right\vert \leq\left\vert \langle
y^{\ast},y-x\rangle\right\vert +\sigma_{U}(y-x)\leq\left\vert \langle y^{\ast
},y-x\rangle\right\vert +\sigma_{U_{0}}(y-x)\leq2\varepsilon,
\]
showing that $z^{\ast}\in\breve{\partial}_{p}^{2\varepsilon}h_{t}(x)$ and so,
\[
y^{\ast}\in z^{\ast}+L^{\bot}+U\subset\breve{\partial}_{p}^{2\varepsilon}%
h_{t}(x)+L^{\bot}+U.
\]
This leads us to
\begin{align*}
\breve{\partial}_{p}^{\varepsilon}(f_{t}+\mathrm{I}_{\overline{L\cap
\operatorname*{dom}f}})(x)  &  \subset\bigcap\nolimits_{U\in\mathcal{N}%
_{X^{\ast}},U\subset U_{0}}\left(  \breve{\partial}_{p}^{2\varepsilon}%
h_{t}(x)+L^{\bot}+U\right) \\
&  =\bigcap\nolimits_{U\in\mathcal{N}_{X^{\ast}}}\left(  \breve{\partial}%
_{p}^{2\varepsilon}h_{t}(x)+L^{\bot}+U\right) \\
&  =\operatorname*{cl}\left(  \breve{\partial}_{p}^{2\varepsilon}%
h_{t}(x)+L^{\bot}\right)  ,
\end{align*}
and (\ref{li}) follows.

Now, by applying (\ref{marco1}), and using (\ref{li}), we obtain
\begin{align*}
\partial f(x)  &  =\bigcap_{_{\substack{\varepsilon>0,\text{ }p\in
\mathcal{P}\\L\in\mathcal{F}(x)}}}\overline{\operatorname*{co}}\left\{
\bigcup\limits_{t\in T_{\varepsilon}(x)}\breve{\partial}_{p}^{\varepsilon
}(f_{t}+\mathrm{I}_{\overline{L\cap\operatorname*{dom}f}})(x)\right\} \\
&  \subset\bigcap_{_{\substack{\varepsilon>0,\text{ }p\in\mathcal{P}%
\\L\in\mathcal{F}(x)}}}\overline{\operatorname*{co}}\left(  \bigcup
\limits_{t\in T_{\varepsilon}(x)}\operatorname*{cl}\left(  \breve{\partial
}_{p}^{2\varepsilon}(f_{t}+\mathrm{I}_{\overline{\operatorname*{dom}f}%
})(x)+L^{\bot}\right)  \right) \\
&  \subset\bigcap_{_{\substack{\varepsilon>0,\text{ }p\in\mathcal{P}%
\\L\in\mathcal{F}(x)}}}\overline{\operatorname*{co}}\left(  \operatorname*{cl}%
\left(  \bigcup\limits_{t\in T_{2\varepsilon}(x)}\breve{\partial}%
_{p}^{2\varepsilon}(f_{t}+\mathrm{I}_{\overline{\operatorname*{dom}f}%
})(x)+L^{\bot}\right)  \right) \\
&  =\bigcap_{_{\substack{\varepsilon>0,\text{ }p\in\mathcal{P}\\L\in
\mathcal{F}(x)}}}\overline{\operatorname*{co}}\left(  \bigcup\limits_{t\in
T_{2\varepsilon}(x)}\breve{\partial}_{p}^{2\varepsilon}(f_{t}+\mathrm{I}%
_{\overline{\operatorname*{dom}f}})(x)+L^{\bot}\right) \\
&  =\bigcap_{_{\varepsilon>0,\text{ }p\in\mathcal{P}}}\overline
{\operatorname*{co}}\left(  \bigcup\limits_{t\in T_{2\varepsilon}(x)}%
\breve{\partial}_{p}^{2\varepsilon}(f_{t}+\mathrm{I}_{\overline
{\operatorname*{dom}f}})(x)\right)  ,
\end{align*}
where in the last equality we used Lemma (\ref{ar}).
\end{dem}

\section{The continuous case}

In this section we get reduced versions of\ the previous formulas, (\ref{FVB})
and (\ref{marco1}), when\ continuity assumptions are imposed. We work in the
general setting of a lcs space $X,$ where $\mathcal{P}$ is the family of
continuous seminorms in $X$. The main result\ comes next.

\begin{theo}
\label{propinter}Let\ $f_{t}\in\Lambda(X),$ $t\in T,$ and $f=\sup_{t\in
T}f_{t}.$ Assume that $f$ is finite and continuous at some point.\ Then for
every $x\in X$
\begin{equation}
\partial f(x)=\mathrm{N}_{\operatorname*{dom}f}(x)+\bigcap_{\varepsilon
>0,\text{ }p\in\mathcal{P}}\overline{\operatorname*{co}}\left\{
\bigcup\limits_{t\in\overline{T}_{\varepsilon}(x)}\breve{\partial}%
_{p}^{\varepsilon}(\operatorname*{cl}f_{t})(x)\right\}  , \label{marco2}%
\end{equation}
where%
\[
\overline{T}_{\varepsilon}(x):=\{t\in T\mid(\operatorname*{cl}f_{t})(x)\geq
f(x)-\varepsilon\}.
\]

\end{theo}

\begin{dem}
We assume without loss of generality\ that $x=\theta\in\operatorname*{dom}f$
and $f(\theta)=0.$ We fix $x^{\ast}\in\partial f(\theta)$ and let $x_{0}%
\in\operatorname*{dom}f$ be such that $f$ is continuous at $x_{0},$ that is,
there exist $W\in\mathcal{N}_{X}$ and $m\in\mathbb{R}_{+}$ such that
\begin{equation}
x_{0}+W\subset\operatorname*{dom}f\text{ \ and }\sup_{t\in T,\text{ }w\in
W}f_{t}(x_{0}+w)\leq m; \label{cw}%
\end{equation}
hence, all the $f_{t}$'s are continuous at $x_{0}$ \cite[Lemma 2.2.8]%
{ZalinescuBook}.

\emph{Step1.} Let us first suppose that all the $f_{t}$'s are lsc, so that
$\overline{T}_{\varepsilon}(x)=T_{\varepsilon}(x).$ Fix $\varepsilon\in\left]
0,1\right[  $ and $p\in\mathcal{P}.$

We introduce the continuous seminorm $\tilde{p}:=\max\{p,p_{W}\},$ where
$p_{W}$ is\ the Minkowski gauge of $W.$ Then, by Theorem \ref{sinl},
\begin{equation}
\partial f(\theta)=\bigcap_{_{\delta>0,\text{ }q\in\mathcal{P}}}%
\overline{\operatorname*{co}}\left\{  \bigcup\limits_{t\in T_{\delta}(\theta
)}\breve{\partial}_{q}^{\delta}(f_{t}+\mathrm{I}_{\overline
{\operatorname*{dom}f}})(\theta)\right\}  \subset\overline{\operatorname*{co}%
}\left\{  \bigcup\limits_{t\in T_{\varepsilon}(\theta)}\breve{\partial
}_{\tilde{p}}^{\varepsilon}(f_{t}+\mathrm{I}_{\overline{\operatorname*{dom}f}%
})(\theta)\right\}  . \label{30}%
\end{equation}

\emph{Step 2.} Choose $V_{0}\in\mathcal{N}_{X^{\ast}}$ satisfying\
\begin{equation}
\sigma_{V_{0}}(x_{0})\leq1, \label{sig}%
\end{equation}
and take $V\in\mathcal{N}_{X^{\ast}}$ such that $V\subset V_{0};$ hence,
$\sigma_{V}(x_{0})\leq\sigma_{V_{0}}(x_{0})\leq1,$ and by (\ref{30})%
\[
x^{\ast}\in\operatorname*{co}\left\{  \bigcup\limits_{t\in T_{\varepsilon
}(\theta)}\breve{\partial}_{\tilde{p}}^{\varepsilon}(f_{t}+\mathrm{I}%
_{\overline{\operatorname*{dom}f}})(\theta)\right\}  +V.
\]
In other words, there exist\ $k\in\mathbb{N},$ $(\lambda_{1},\cdots
,\lambda_{k})\in\Delta_{k}$, $t_{1},\cdots,t_{k}\in T_{\varepsilon}(\theta),$
and $z_{i}^{\ast}\in\breve{\partial}_{\tilde{p}}^{\varepsilon}(f_{t_{i}%
}+\mathrm{I}_{\overline{\operatorname*{dom}f}})(\theta),$ $i\in\overline
{1,k},$ such that
\[
x^{\ast}\in\lambda_{1}z_{1}^{\ast}+\cdots+\lambda_{k}z_{k}^{\ast}+V.
\]
Moreover, taking into account Moreau-Rockafellar sum rule (recall that the
$f_{t}$'s are continuous at $x_{0}\in\operatorname*{dom}f$), and using the
definition of $\breve{\partial}_{\tilde{p}}^{\varepsilon},$ for each
$i\in\overline{1,k}$ there exist $y_{i}\in\overline{\operatorname*{dom}f}$,
with $\tilde{p}(y_{i})\leq\varepsilon,$ and elements $y_{i}^{\ast}\in\partial
f_{t_{i}}(y_{i})$ and $v_{i}^{\ast}\in\mathrm{N}_{\overline
{\operatorname*{dom}f}}(y_{i})$ such that $z_{i}^{\ast}=y_{i}^{\ast}%
+v_{i}^{\ast},$
\begin{equation}
\left\vert \left\langle y_{i}^{\ast}+v_{i}^{\ast},y_{i}\right\rangle
\right\vert \leq\varepsilon, \label{er1}%
\end{equation}
and
\begin{equation}
\left\vert f_{t_{i}}(y_{i})-f_{t_{i}}(\theta)\right\vert =\left\vert
(f_{t_{i}}+\mathrm{I}_{\overline{\operatorname*{dom}f}})(y_{i})-(f_{t_{i}%
}+\mathrm{I}_{\overline{\operatorname*{dom}f}})(\theta)\right\vert
\leq\varepsilon; \label{er}%
\end{equation}
consequently, $y_{i}$ and $y_{i}^{\ast}$ satisfy the following\ relations
\begin{equation}
y_{i}^{\ast}\in\partial_{3\varepsilon}f(\theta)\ \text{and }\left\vert
\left\langle y_{i}^{\ast},y_{i}\right\rangle \right\vert \leq\varepsilon,
\label{te}%
\end{equation}
and
\begin{equation}
y_{i}^{\ast}\in\breve{\partial}_{\tilde{p}}^{\varepsilon}f_{t_{i}}(\theta).
\label{gg}%
\end{equation}
Indeed, for any $u\in\operatorname*{dom}f$ the information we have above on
$y_{i}$ and $y_{i}^{\ast}$ leads, on\ one hand,\ to
\begin{align*}
\left\langle y_{i}^{\ast},u\right\rangle  &  =\left\langle y_{i}^{\ast
},u-y_{i}\right\rangle +\left\langle y_{i}^{\ast}+v_{i}^{\ast},y_{i}%
\right\rangle +\left\langle v_{i}^{\ast},-y_{i}\right\rangle \\
&  \leq f_{t_{i}}(u)-f_{t_{i}}(y_{i})+\varepsilon\\
&  \leq f(u)-f_{t_{i}}(\theta)+2\varepsilon\\
&  \leq f(u)-f(\theta)+3\varepsilon,
\end{align*}
which shows that $y_{i}^{\ast}\in\partial_{3\varepsilon}f(\theta).$ The second
inequality in (\ref{te}) also follows since
\[
\left\langle y_{i}^{\ast},-y_{i}\right\rangle \leq f_{t_{i}}(\theta)-f_{t_{i}%
}(y_{i})\leq\varepsilon,
\]
and, using (\ref{er1}),
\[
\left\langle y_{i}^{\ast},y_{i}\right\rangle =\left\langle y_{i}^{\ast}%
+v_{i}^{\ast},y_{i}\right\rangle +\left\langle v_{i}^{\ast},-y_{i}%
\right\rangle \leq\varepsilon.
\]
This shows that $\left\vert \left\langle y_{i}^{\ast},y_{i}\right\rangle
\right\vert \leq\varepsilon,$ which together with (\ref{er}), and the
relations $\tilde{p}(y_{i})\leq\varepsilon$ and $y_{i}^{\ast}\in\partial
f_{t_{i}}(y_{i}),$ lead us to (\ref{gg}).

We also have
\begin{equation}
v_{i}^{\ast}\in\mathrm{N}_{\overline{\operatorname*{dom}f}}^{2\varepsilon
}(\theta). \label{teb}%
\end{equation}
Actually, for any $u\in\operatorname*{dom}f,$ by (\ref{er1}) and (\ref{te}) we
have
\begin{align*}
\langle v_{i}^{\ast},u\rangle &  =\langle v_{i}^{\ast},u-y_{i}\rangle+\langle
v_{i}^{\ast},y_{i}\rangle\leq\langle v_{i}^{\ast},y_{i}\rangle\\
&  =\langle v_{i}^{\ast}+y_{i}^{\ast},y_{i}\rangle-\langle y_{i}^{\ast}%
,y_{i}\rangle\leq2\varepsilon,
\end{align*}
and this yields (\ref{teb}).

Also, since $p_{W}(y_{i})\leq\tilde{p}(y_{i})\leq\varepsilon,$ we infer
that\ $y_{i}\in\varepsilon W$ and so, using (\ref{cw}),%
\[
x_{0}+y_{i}\in x_{0}+\varepsilon W\subset x_{0}+W\subset\operatorname*{dom}f,
\]
that is, since $v_{i}^{\ast}\in\mathrm{N}_{\overline{\operatorname*{dom}f}%
}(y_{i}),$\
\begin{equation}
\left\langle v_{i}^{\ast},x_{0}\right\rangle =\left\langle v_{i}^{\ast}%
,(x_{0}+y_{i})-y_{i}\right\rangle \leq0. \label{cw2}%
\end{equation}

\emph{Step 3.} Now, we denote
\[
y_{V}^{\ast}:=\lambda_{1}y_{1}^{\ast}+\cdots+\lambda_{k}y_{k}^{\ast},\text{
}v_{V}^{\ast}:=\lambda_{1}v_{1}^{\ast}+\cdots+\lambda_{k}v_{k}^{\ast};
\]
hence, by (\ref{te}), (\ref{cw2}) and (\ref{gg}),
\begin{equation}
y_{V}^{\ast}\in\partial_{3\varepsilon}f(\theta)\cap\operatorname*{co}\left\{
\bigcup\limits_{t\in T_{\varepsilon}(\theta)}\breve{\partial}_{\tilde{p}%
}^{\varepsilon}f_{t}(\theta)\right\}  , \label{gg1}%
\end{equation}
and
\begin{equation}
v_{V}^{\ast}\in\mathrm{N}_{\overline{\operatorname*{dom}f}}^{2\varepsilon
}(\theta),\text{ }\left\langle v_{V}^{\ast},x_{0}\right\rangle \leq0,
\label{gg2}%
\end{equation}
at the same time as\
\begin{equation}
v_{V}^{\ast}\in x^{\ast}-y_{V}^{\ast}+V, \label{31}%
\end{equation}
entailing that, due to (\ref{sig}),%
\begin{align}
\langle y_{V}^{\ast},x_{0}\rangle &  =\langle y_{V}^{\ast}+v_{V}^{\ast}%
,x_{0}\rangle-\langle v_{V}^{\ast},x_{0}\rangle\nonumber\\
&  \geq\langle x^{\ast},x_{0}\rangle-\sigma_{V}(x_{0})\nonumber\\
&  \geq\langle x^{\ast},x_{0}\rangle-\sigma_{V_{0}}(x_{0})\geq\langle x^{\ast
},x_{0}\rangle-1. \label{32}%
\end{align}
Consequently, since $y_{V}^{\ast}\in\partial_{3\varepsilon}f(\theta),$ we
obtain\ that, for all $w\in W,$
\begin{equation}
\langle y_{V}^{\ast},x_{0}+w\rangle\leq f(x_{0}+w)+3\varepsilon\leq
m+3\varepsilon, \label{pol1}%
\end{equation}
and so, because of\ (\ref{32}), we find\ some $r>0$ such that
\begin{equation}
\langle y_{V}^{\ast},w\rangle\leq-\langle y_{V}^{\ast},x_{0}\rangle
+m+3\varepsilon\leq-\langle x^{\ast},x_{0}\rangle+1+m+3\varepsilon\leq r;
\label{ac}%
\end{equation}
that is, considering the natural partial order in $\mathcal{N}_{X^{\ast}},$ by
Alaoglu-Bourbaki theorem the net $(y_{V}^{\ast})_{V\in\mathcal{N}_{X^{\ast}%
},\text{ }V\subset V_{0}}$ is contained in the weak*-compact set $rW^{\circ}$.
Hence, we may assume that $(y_{V}^{\ast})_{V}$ weak*-converges to some element
$y_{\varepsilon}^{\ast}\in X^{\ast}$ such that (recall (\ref{gg1}))
\begin{equation}
y_{\varepsilon}^{\ast}\in rW^{\circ}\cap\overline{\operatorname*{co}}\left\{
\bigcup\limits_{t\in T_{\varepsilon}(\theta)}\breve{\partial}_{\tilde{p}%
}^{\varepsilon}f_{t}(\theta)\right\}  \subset rW^{\circ}\cap\overline
{\operatorname*{co}}\left\{  \bigcup\limits_{t\in T_{\varepsilon}(x)}%
\breve{\partial}_{\tilde{p}}^{\varepsilon}f_{t}(\theta)\right\}  . \label{gg3}%
\end{equation}
It also follows that the corresponding (sub)net\ $(v_{V}^{\ast})_{V\in
\mathcal{N}_{X^{\ast}},\text{ }V\subset V_{0}}$ weak*-converges to the element
$v_{\varepsilon}^{\ast}:=x^{\ast}-y_{\varepsilon}^{\ast}$ and that
$v_{\varepsilon}^{\ast}\in\mathrm{N}_{\overline{\operatorname*{dom}f}%
}^{2\varepsilon}(\theta).$ Indeed, given $U_{1}\in\mathcal{N}_{X^{\ast}},$ we
choose $U_{0}\in\mathcal{N}_{X^{\ast}}$ such that $U_{0}+U_{0}\subset U_{1}.$
Since $U:=v_{\varepsilon}^{\ast}+U_{0}\in\mathcal{N}_{v_{\varepsilon}^{\ast}%
},$ we have $x^{\ast}-U\in\mathcal{N}_{y_{\varepsilon}^{\ast}}$ and so, there
is some $V_{1}\in\mathcal{N}_{X^{\ast}}$ such that $y_{V^{\prime}}^{\ast}\in
x^{\ast}-U$ for all $V^{\prime}\subset V_{1}$; that is, $\theta\in
y_{V^{\prime}}^{\ast}-x^{\ast}+U$ and, using (\ref{31}),
\[
v_{V^{\prime}}^{\ast}\in v_{V^{\prime}}^{\ast}+y_{V^{\prime}}^{\ast}-x^{\ast
}+U\subset V^{\prime}+U.
\]
Then we get,\ since $V^{\prime}\cap U_{0}\subset V^{\prime}$ and
$U=v_{\varepsilon}^{\ast}+U_{0},$
\[
v_{V^{\prime}}^{\ast}\in V^{\prime}\cap U_{0}+U=v_{\varepsilon}^{\ast
}+V^{\prime}\cap U_{0}+U_{0}\subset v_{\varepsilon}^{\ast}+U_{1},
\]
showing that $(v_{V}^{\ast})_{V\in\mathcal{N}_{X^{\ast}},\text{ }V\subset
V_{0}}$ weak*-converges to $v_{\varepsilon}^{\ast}\in\mathrm{N}_{\overline
{\operatorname*{dom}f}}^{2\varepsilon}(\theta),$ by (\ref{gg2}).

\emph{Step 4. }Again, since $(y_{\varepsilon}^{\ast})_{\varepsilon}\subset
rW^{\circ}$ we may assume that the net $(y_{\varepsilon}^{\ast})_{\varepsilon
}$ weak*-converge to some $y^{\ast},$ which by (\ref{gg3}) satisfies
\[
y^{\ast}\in\bigcap\limits_{\varepsilon>0}\overline{\operatorname*{co}}\left\{
\bigcup\limits_{t\in T_{\varepsilon}(\theta)}\breve{\partial}_{\tilde{p}%
}^{\varepsilon}f_{t}(\theta)\right\}  .
\]
Moreover, since $p$ was arbitrarily fixed in $\mathcal{P}$, we deduce that
\[
y^{\ast}\in\bigcap\limits_{\varepsilon>0,\text{ }p\in\mathcal{P}}%
\overline{\operatorname*{co}}\left\{  \bigcup\limits_{t\in T_{\varepsilon
}(\theta)}\breve{\partial}_{\tilde{p}}^{\varepsilon}f_{t}(\theta)\right\}  .
\]
In addition, since $v_{\varepsilon}^{\ast}=x^{\ast}-y_{\varepsilon}^{\ast},$
by (\ref{gg2}) it follows that $(v_{\varepsilon}^{\ast})_{\varepsilon}$ also
weak*-converges to
\[
v^{\ast}=x^{\ast}-y^{\ast}\in\mathrm{N}_{\overline{\operatorname*{dom}f}%
}(\theta)=\mathrm{N}_{\operatorname*{dom}f}(\theta).
\]
Consequently,
\[
x^{\ast}=v^{\ast}+y^{\ast}\in\mathrm{N}_{\operatorname*{dom}f}(\theta
)+\bigcap\limits_{\varepsilon>0,\text{ }p\in\mathcal{P}}\overline
{\operatorname*{co}}\left\{  \bigcup\limits_{t\in T_{\varepsilon}(\theta
)}\breve{\partial}_{\tilde{p}}^{\varepsilon}f_{t}(\theta)\right\}  ,
\]
yielding\ the inclusion \textquotedblleft$\subset$\textquotedblright.

\emph{Step 5.} We consider the case when the $f_{t}$'s are not necessarily
lsc. Due to the continuity assumption on $f,$ entailing the continuity of the
$f_{t}$'s at $x_{0},$ according to \cite[Corollary 9(ii)]{HanLopZal2008} we
have
\[
\operatorname*{cl}f=\sup_{t\in T}\operatorname*{cl}f_{t}.
\]
If $x\in\operatorname*{dom}f$ is such that $\partial f(x)\neq\emptyset,$ then
$\partial f(x)=\partial(\operatorname*{cl}f)(x)$ and $f(x)=(\operatorname*{cl}%
f)(x).$ Thus, from Step 4 applied to the $(\operatorname*{cl}f_{t})$'s, we
obtain
\[
\partial f(x)=\partial(\operatorname*{cl}f)(x)=\mathrm{N}_{\operatorname*{dom}%
(\operatorname*{cl}f)}(x)+\bigcap\limits_{\varepsilon>0,\text{ }%
p\in\mathcal{P}}\overline{\operatorname*{co}}\left\{  \bigcup\limits_{t\in
\overline{T}_{\varepsilon}(x)}\breve{\partial}_{\tilde{p}}^{\varepsilon
}(\operatorname*{cl}f_{t})(x)\right\}  ,
\]
and the inclusion \textquotedblleft$\subset$\textquotedblright\ follows since
$\mathrm{N}_{\operatorname*{dom}(\operatorname*{cl}f)}(x)\subset
\mathrm{N}_{\operatorname*{dom}f}(x).$

\emph{Step 6.} We show the opposite inclusion "$\supset$". Given
$x\in\operatorname*{dom}f$ $(\subset\operatorname*{dom}f_{t}),$ by (\ref{cco})
we have $\breve{\partial}_{\tilde{p}}^{\varepsilon}\!(\operatorname*{cl}%
f_{t})(x)\subset\partial_{2\varepsilon}(\operatorname*{cl}f_{t})(x).$ Then,
since it can be easily proved that\ $\partial_{2\varepsilon}%
(\operatorname*{cl}f_{t})(x)\subset\partial_{3\varepsilon}f(x)$ for all
$t\in\overline{T}_{\varepsilon}(x),$ we deduce that for all $\varepsilon>0$
and $p\in\mathcal{P}$
\[
\overline{\operatorname*{co}}\left\{  \bigcup\limits_{t\in\overline
{T}_{\varepsilon}(x)}\breve{\partial}_{\tilde{p}}^{\varepsilon}%
(\operatorname*{cl}f_{t})(x)\right\}  \subset\partial_{3\varepsilon}f(x).
\]
Therefore
\begin{align*}
\mathrm{N}_{\operatorname*{dom}f}(x)+\bigcap\limits_{\varepsilon>0,\text{
}p\in\mathcal{P}}\overline{\operatorname*{co}}\left\{  \bigcup\limits_{t\in
\overline{T}_{\varepsilon}(x)}\breve{\partial}_{\tilde{p}}^{\varepsilon
}(\operatorname*{cl}f_{t})(x)\right\}   &  \subset\mathrm{N}%
_{\operatorname*{dom}f}(x)+\bigcap_{\varepsilon>0}\partial_{\varepsilon}f(x)\\
&  =\mathrm{N}_{\operatorname*{dom}f}(x)+\partial f(x)=\partial f(x).
\end{align*}

\end{dem}

\begin{rem}
\emph{As can be observed from the proof of Theorem \ref{propinter}, the
continuity assumption on }$f$\emph{ used there can be weakened to the
existence of small }$\varepsilon>0$\emph{ such that the function }$\sup
_{t\in\overline{T}_{\varepsilon}(x)}f_{t}$\emph{ is finite and continuous at
some point.}
\end{rem}

The following corollary is straightforward from Theorem \ref{propinter}.

\begin{cor}
Let\ $f_{t}\in\Gamma_{0}(X),$ $t\in T,$ and $f=\sup_{t\in T}f_{t}.$ Assume
that $f$ is finite and continuous at some point.\ Then for every $x\in X$
\[
\partial f(x)=\mathrm{N}_{\operatorname*{dom}f}(x)+\bigcap_{\varepsilon
>0,\text{ }p\in\mathcal{P}}\overline{\operatorname*{co}}\left\{
\bigcup\limits_{t\in T_{\varepsilon}(x)}\breve{\partial}_{\tilde{p}%
}^{\varepsilon}f_{t}(x)\right\}  .
\]

\end{cor}

We obtain the classical result of Valadier \cite[Theorem 1]{Valadier}:

\begin{cor}
Let\ $f_{t}\in\Lambda(X),$ $t\in T,$ and $f=\sup_{t\in T}f_{t}.$ Assume that
$f$ is continuous at $x\in\operatorname*{dom}f.$\ Then
\[
\partial f(x)=\bigcap_{\varepsilon>0,\text{ }p\in\mathcal{P}}\overline
{\operatorname*{co}}\left\{  \bigcup\limits_{t\in T_{\varepsilon}(x),\text{
}p(y-x)\leq\varepsilon}\partial f_{t}(y)\right\}  .
\]

\end{cor}

\begin{dem}
We assume again that $x=\theta\in\operatorname*{dom}f$ and $f(\theta)=0.$
Since $f$ is continuous at $\theta,$ by Theorem \ref{propinter} we have
\begin{align}
\partial f(\theta)  &  =\mathrm{N}_{\operatorname*{dom}f}(\theta
)+\bigcap_{\varepsilon>0,\text{ }p\in\mathcal{P}}\overline{\operatorname*{co}%
}\left\{  \bigcup\limits_{t\in\overline{T}_{\varepsilon}(\theta)}%
\breve{\partial}_{\tilde{p}}^{\varepsilon}(\operatorname*{cl}f_{t}%
)(\theta)\right\} \nonumber\\
&  =\bigcap_{\varepsilon>0,\text{ }p\in\mathcal{P}}\overline
{\operatorname*{co}}\left\{  \bigcup\limits_{t\in\overline{T}_{\varepsilon
}(\theta)}\breve{\partial}_{\tilde{p}}^{\varepsilon}(\operatorname*{cl}%
f_{t})(\theta)\right\}  . \label{u}%
\end{align}
On\ one hand, the continuity of $f$ at $\theta$ implies the continuity of all
the $f_{t}$'s at $\theta.$ So,
\[
(\operatorname*{cl}f_{t})(\theta)=f_{t}(\theta)\text{ for all }t\in T,
\]
and we get
\begin{equation}
\overline{T}_{\varepsilon}(\theta)=T_{\varepsilon}(\theta). \label{v}%
\end{equation}
On the other hand, there exist $U\in\mathcal{N}_{X}$ and $m\geq0$ such that,
for all $u\in U,$
\begin{equation}
f_{t}(u)\leq f(u)\leq m, \label{ct}%
\end{equation}
and so, by \cite[Lemma 2.2.8]{ZalinescuBook},
\begin{equation}
\left\vert f_{t}(u)-f_{t}(\theta)\right\vert \leq mp_{U}(u). \label{ct2}%
\end{equation}
This entails the continuity of the $f_{t}$'s at\ the points of $U.$ Hence,
(\ref{u}) together with (\ref{v}) give us, for all $p\in\mathcal{P}$,
\[
\partial f(\theta)\subset\bigcap_{\varepsilon>0}\overline{\operatorname*{co}%
}\left\{  \bigcup\limits_{t\in T_{\varepsilon}(\theta)}\breve{\partial
}_{\tilde{p}}^{\varepsilon}(\operatorname*{cl}f_{t})(\theta)\right\}  ,
\]
where $\tilde{p}:=\max\{p,p_{U}\}.$ Observe that if $y^{\ast}\in
\breve{\partial}_{\tilde{p}}^{\varepsilon}(\operatorname*{cl}f_{t})(\theta)$,
for $t\in T_{\varepsilon}(\theta),$ then there exists $y\in X$ satisfying
$\tilde{p}(y)\leq\varepsilon$ such that $y^{\ast}\in\partial
(\operatorname*{cl}f_{t})(y).$ But $p_{U}(y)\leq\varepsilon$ ensures that
$y\in\varepsilon U\subset U$, and so, by (\ref{ct}), $f_{t}$ is continuous at
$y$ so that $y^{\ast}\in\partial(\operatorname*{cl}f_{t})(y)=\partial
f_{t}(y).$ This gives rise to\
\[
\partial f(\theta)\subset\bigcap_{\varepsilon>0}\overline{\operatorname*{co}%
}\left\{  \bigcup\limits_{t\in T_{\varepsilon}(\theta),\text{ }p(y)\leq
\varepsilon}\partial f_{t}(y)\right\}  ,
\]
and the inclusion \textquotedblleft$\subset$\textquotedblright\ follows by
intersecting over $p\in\mathcal{P}.$

To show the converse inclusion "$\supset$", we take
\[
x^{\ast}\in\bigcap_{\varepsilon>0,\text{ }p\in\mathcal{P}}\overline
{\operatorname*{co}}\left\{  \bigcup\limits_{t\in T_{\varepsilon}%
(\theta),\text{ }p(y)\leq\varepsilon}\partial f_{t}(y)\right\}  .
\]
Then, given $\varepsilon\in\left]  0,\frac{1}{2}\right[  $ and $p\in
\mathcal{P}$, we get
\begin{equation}
x^{\ast}\in\overline{\operatorname*{co}}\left\{  \bigcup\limits_{t\in
T_{\varepsilon}(\theta),\text{ }\tilde{p}(y)\leq\varepsilon}\partial
f_{t}(y)\right\}  , \label{vo}%
\end{equation}
where $\tilde{p}:=\max\{p,p_{U}\}$ $(\in\mathcal{P})$, with $U$ being defined
as in (\ref{ct}). Observe that if\ $y^{\ast}\in\partial f_{t}(y)$ for some
$t\in T_{\varepsilon}(\theta)$ and $y\in X$ such that $\tilde{p}%
(y)\leq\varepsilon,$ then $y\in\varepsilon U\subset\frac{1}{2}U$ and, taking
into account (\ref{ct2})$,$%
\[
\left\vert f_{t}(y)-f_{t}(\theta)\right\vert \leq mp_{U}(y)\leq m\tilde
{p}(y)\leq m\varepsilon,
\]%
\[
\left\vert f_{t}(2y)-f_{t}(\theta)\right\vert \leq mp_{U}(2y)\leq2m\tilde
{p}(y)\leq2m\varepsilon.
\]
Consequently, for every $z\in\operatorname*{dom}f,$
\begin{align*}
\left\langle y^{\ast},z\right\rangle  &  =\left\langle y^{\ast}%
,z-y\right\rangle +\left\langle y^{\ast},y\right\rangle \\
&  =\left\langle y^{\ast},z-y\right\rangle +\left\langle y^{\ast
},2y-y\right\rangle \\
&  \leq f_{t}(z)-f_{t}(y)+f_{t}(2y)-f_{t}(y)\\
&  =(f_{t}(z)-f_{t}(\theta))+(f_{t}(\theta)-f_{t}(y))+(f_{t}(2y)-f_{t}%
(\theta))+(f_{t}(\theta)-f_{t}(y))\\
&  \leq f_{t}(z)-f_{t}(\theta)+4m\varepsilon,
\end{align*}
and we get, since $t\in T_{\varepsilon}(\theta),$%
\[
\left\langle y^{\ast},z\right\rangle \leq f(z)-f(\theta)+4m\varepsilon
+\varepsilon,
\]
so that $y^{\ast}\in\partial_{(4m+1)\varepsilon}f(\theta).$ Thus, from
(\ref{vo}) it follows that
\[
x^{\ast}\in\bigcap_{0<\varepsilon<\frac{1}{2}}\partial_{(4m+1)\varepsilon
}f(\theta)=\partial f(\theta),
\]
and the proof is complete.
\end{dem}

\end{document}